\documentclass[%
a4paper
]%
{article}

\usepackage{amsmath}
\usepackage{amssymb}
\usepackage{amsthm}
\usepackage{amscd}
\usepackage{xspace}

\newtheoremstyle{thm}{3pt}{3pt}{\itshape}{}{\bfseries}{}{.5em}{}
 \newtheorem{theorem}{Theorem}[section]
 \newtheorem{lemma}[theorem]{Lemma}
 \newtheorem{proposition}[theorem]{Proposition}
 \newtheorem{defn}[theorem]{Definition}
 \newtheorem{corollary}[theorem]{Corollary}

\newtheoremstyle{thmsub}{3pt}{3pt}{\upshape}{}{\bfseries}{}{.5em}{}
 \newtheorem{remark}[theorem]{Remark}

\newcommand{\diff}[2]{\ensuremath\frac{d#1}{d#2}}
\newcommand{\C}{\ensuremath{\mathbb{C}}\xspace}
\newcommand{\CP}{\ensuremath{\mathbb{CP}}\xspace}
\newcommand{\F}{\ensuremath{\mathbb{F}}\xspace}
\newcommand{\R}{\ensuremath{\mathbb{R}}\xspace}
\newcommand{\Z}{\ensuremath{\mathbb{Z}}\xspace}
\newcommand{\T}{\ensuremath{\mathbb{T}}\xspace}
\newcommand{\Ci}{\ensuremath{C^\infty}\xspace}
\newcommand{\restrict}{\!\!\mid}
\newcommand{\conj}[1]{\overline{#1}}
\newcommand{\exterior}{\Lambda}
\newcommand{\m}[1]{\ensuremath{\mathcal{#1}}\xspace}
\newcommand{\mf}[1]{\ensuremath{\mathfrak{#1}}\xspace}
\newcommand{\ssetminus}{\ensuremath{\!\smallsetminus\!}}
\newcommand{\abs}[1]{\left\lvert#1\right\rvert}
\newcommand{\norm}[1][\cdot]{\left\lVert#1\right\rVert}
\newcommand{\ip}[2]{\langle#1,#2\rangle}
\newcommand{\nmref}[1]{\ref{#1}}
\newcommand{\nmeqref}[1]{\eqref{#1}}
\newcommand{\ipv}[2]{\left( #1, #2 \right)}
\newcommand{\ipvc}{\ipv{\cdot}{\cdot}}
\newcommand{\ipc}{\ip{\cdot}{\cdot}}
\newcommand{\bbh}{\ensuremath{\mathbb{H}}}
\newcommand{\dirac}{\partial\hspace*{-1.2ex}/}
\newcommand{\pol}{\text{pol}}
\newcommand{\wotimes}{\widetilde{\otimes}}
\newcommand{\ad}{\text{ad}}
\newcommand{\per}{\text{per}}
\newcommand{\type}{\text{abc}}
\newcommand{\res}{\text{res}}

\newcommand{\ttype}{\emph{abc}\xspace}
\newcommand{\tper}{\emph{per}\xspace}
\newcommand{\tpol}{\emph{pol}\xspace}

\DeclareMathOperator{\Diff}{Diff}
\DeclareMathOperator{\pin}{Pin}
\DeclareMathOperator{\Ad}{Ad}
\DeclareMathOperator{\gl}{Gl}
\DeclareMathOperator{\ind}{Index}
\DeclareMathOperator{\tr}{Tr}
\DeclareMathOperator{\spin}{Spin}
\DeclareMathOperator{\gr}{Gr}

\numberwithin{equation}{section}

\allowdisplaybreaks

\begin{document}

\title{The Geometry of the Loop Space and a Construction of a Dirac
  Operator}
\date{\today}
\author{Andrew Stacey}

\maketitle

\begin{abstract}
We describe a construction of fibrewise inner products on the
cotangent bundle of the smooth free loop space of a Riemannian
manifold.  Using this inner product, we construct an operator over the
loop space of a string manifold which is directly analogous to the
Dirac operator of a spin manifold.
\end{abstract}

\tableofcontents

\newpage
\section{Introduction}
\label{sec:intro}

The problem addressed in this paper is that of constructing an inner
product on the cotangent bundle of the manifold of smooth unbased
loops in a smooth finite dimensional Riemannian manifold.  The main
motivation for this is the problem of constructing for the loop space
an analogue of the Dirac operator of a finite dimensional spin
manifold.

In this introduction we start with an overview of the construction of
the Dirac operator in finite dimensions, explain what can and cannot
be generalised to infinite dimensions, and show how an inner product
on the cotangent bundle solves the problems that occur.  We follow
this with a short discussion on the connection between inner products
and Hilbert completions and explain what exactly we aim to construct
in the paper.  The main part of this introduction finishes with an
outline of the method of construction.

\subsection{The Dirac Operator in Finite Dimensions}
\label{sec:introfinite}

The construction of the Dirac operator is the main motivation for the
construction of the inner product on the cotangent bundle of the loop
space so we explain this first.  We start with the construction in
finite dimensions.  As this is laid out in detail elsewhere we shall
focus on the pieces that lead to the difficulties in infinite
dimensions.  For more on the details of the construction in finite dimensions see~\cite{hlmm}.
For details of the spin representation in all dimensions
see~\cite{rppr}.

There are two methods of constructing the Dirac operator in finite
dimensions.  Both follow the same general outline and have the same
initial data, namely a spin manifold \(M\).  Part of what is meant by
the statement that ``\(M\) is spin'' is that \(M\) is a Riemannian
manifold and so there is an inner product on the tangent bundle.
Since the Riemannian structure defines an isomorphism of the tangent
and cotangent bundles we can transfer this inner product to the
cotangent bundle.

The first step is to construct two finite dimensional unitary vector
bundles over \(M\) called the \emph{spinor bundles} of \(M\).  What is
relevant for our purposes is that the construction starts from a
vector bundle with an inner product.  This is where the two methods
diverge: one starts with the tangent bundle, the other with the
cotangent bundle.  We shall write \(S^+_T\) and \(S^-_T\) for the
bundles constructed from the tangent bundle and \(S^+_{T^*}\),
\(S^-_{T^*}\) for those from the cotangent bundle.  When we wish to
refer to something that holds for both methods, we shall use the
notation \(S^+\) and \(S^-\) for the spinor bundles and \(T^? M\) for
the correct choice of tangent or cotangent bundle.

The key properties of these spinor bundles are the following: first,
there is an operation called \emph{Clifford multiplication} which is a
vector bundle map:
\[
c : T M \otimes S^+_T \to S^-_T, \qquad c : T^* M \otimes S^+_{T^*} \to
S^-_{T^*}.
\]
Second, there is a natural covariant differential operator \(\nabla\)
on \(S^+\) arising from the Levi-Civita connection on \(M\).

The Clifford multiplication map extends in the natural way to a linear
map on sections.  Together with the differential operator, we
therefore have maps:
\begin{equation}
\label{diag:constr}
\begin{CD}
  \Gamma(S^+) @>\nabla>> \Gamma(\m{L}(T M, S^+)) \\
\\
  && \Gamma(T^? M \otimes S^+) @>c>> \Gamma(S^-).
\end{CD}
\end{equation}
Here \(\m{L}(T M, S^+)\) means the bundle with fibres the linear maps
between the corresponding fibres of \(T M\) and \(S^+\).  Note that
the \(T M\) that appears here is definitely \(T M\) and not \(T^? M\).

We wish to compose these maps to define the Dirac operator.  In order
to do so we must find a vertical map to fill the gap.  This is not
difficult.  First, we observe that for finite dimensional spaces \(V\)
and \(W\) the space of linear maps from \(V\) to \(W\),
\(\m{L}(V,W)\), is naturally isomorphic to \(V^* \otimes W\) where
\(V^*\) is the linear dual of \(V\).  Therefore \(\m{L}(T M, S^+)
\cong T^* M \otimes S^+\).  If we are using the cotangent method, we
can stop here as this is the domain of the Clifford multiplication
map.  If we are using the tangent method we must use the inner product
on the tangent bundle to identify \(T M\) with \(T^* M\) and thus \(T
M \otimes S^+\) with \(T^* M \otimes S^+\).  Thus we obtain the Dirac
operator \(\dirac : \Gamma(S^+) \to \Gamma(S^-)\) as one of the
compositions:
\[
\begin{CD}
\text{(cotangent:) } @. \Gamma(S^+_{T^*}) @>\nabla>> \Gamma(\m{L}(T
M, S^+_{T^*})) \\
&&&& @V\cong VV \\
&&&& \Gamma(T^* M \otimes S^+_{T^*}) @>c>> \Gamma(S^-_{T^*}).
\end{CD}
\]

\[
\begin{CD}
\text{(tangent:) } @. \Gamma(S^+_T) @>\nabla>> \Gamma(\m{L}(T
M, S^+_{T})) \\
&&&& @V\cong VV \\
&&&& \Gamma(T^* M \otimes S^+_T) \\
&&&& @V\cong VV \\
&&&& \Gamma(T M \otimes S^+_{T}) @>c>> \Gamma(S^-_{T}).
\end{CD}
\]

The identification of the tangent and cotangent bundles via the inner
product on the tangent bundle defines an isomorphism of the spinor
bundles and thus the two methods lead to isomorphic operators.

At first sight it appears that the tangent method uses only the inner
product on the tangent bundle whilst the cotangent method uses only
the inner product on the cotangent bundle (ignoring, for the moment,
the fact that the inner product on the cotangent bundle was defined
using the one on the tangent bundle).  The first part of that
statement is not strictly true.  The inner product on the cotangent
bundle appears surreptitiously in the identification of the tangent and
cotangent bundles.  The inner product on the tangent bundle defines an
injective map \(T M \to T^* M\) which, for dimension reasons, is an
isomorphism.  In the construction of the Dirac operator it is not this
map which is used but rather its inverse, \(T^* M \to T M\).  Whilst
we can think of this as merely the inverse to the above map, it is
useful to think of it as the natural map coming from the inner product
on the cotangent bundle, \(T^* M \to {T^*}^* M\).

Therefore once the inner product has been transferred to the cotangent
bundle, the cotangent method only uses that inner product while the
tangent method uses both.

\subsection{Generalising to Loop Spaces}
\label{sec:introloop}

We now consider what is known to generalise -- prior to this paper --
from the finite dimensional construction of the Dirac operator to the
case of loop spaces.  Essentially, everything generalises for the
tangent method up to and including diagram~\eqref{diag:constr}.  Thus
for a loop space which is spin there are bundles \(S^+_T, S^-_T \to L
M\) -- now unitary Hilbert bundles -- together with a covariant
differential operator and a Clifford multiplication map as before
(although the differential operator is not as natural as the finite
dimensional one).  These bundles are constructed from the tangent
bundle with its natural inner product.

The cotangent method is dead in the water as it requires an inner
product as part of its initial data and -- prior to this paper -- such
has not been defined.

The next step for the tangent method is to fill in the gap in the
analogous diagram to~\eqref{diag:constr}.  This gap in finite
dimensions was filled in by two maps.  The first of these came from
the natural isomorphism, in finite dimensions, of \(\m{L}(V,W)\) with
\(V^* \otimes W\).  The natural map is \(V^* \otimes W \to
\m{L}(V,W)\), \(f \otimes w \to (v \to f(v) w)\), and this map exists
for any vector spaces.  It is not generally an isomorphism in infinite
dimensions.

In the case that we are dealing with, \(V\) is a \emph{complete,
nuclear, reflexive} space and \(W\) is a Hilbert space.  It is a
remarkable fact that for such spaces the completion of \(V^* \otimes
W\) with respect to the projective tensor product topology \emph{is}
isomorphic to \(\m{L}(V,W)\) under the natural map above.
Essentially, this is because the completion of \(V^* \otimes W\) is
the space of all \emph{compact} maps from \(V\) to \(W\) and under the
assumptions on \(V\) and \(W\), all continuous maps are compact.  We
prove this remarkable isomorphism in proposition~\ref{prop:cmext}.

The Clifford multiplication map extends over the corresponding
completion so we can complete all tensor products in this diagram with
respect to the projective topology.  Thus we can fill in half the gap.
Here, however, the method stalls.  The inner product on the tangent
bundle defines, as before, an injective map \(T L M \to T^* L M\) but
this is not -- and cannot be made to be -- an isomorphism.  This is
purely a linear question and is due to the fact that the model spaces
for the fibres are not isomorphic.  Since we want the inverse of this
map, our construction of the Dirac operator by the tangent method
falls here.

Both methods fail due to the same problem: a lack of an inner product
on the cotangent bundle.  If we had such an object, we could resurrect
the cotangent method -- providing some technical conditions are
satisfied.  Whereupon  the gap in~\eqref{diag:constr} for
the cotangent method can be filled by the remarkable isomorphism;
ergo: the cotangent method will yield a Dirac operator.  We could also
restart the stalled tangent method since the inner product on the
cotangent bundle would define the injective map \(T^* L M \to T L M\)
(the model space is reflexive so \({T^*}^* L M \cong T L M\) just as
in finite dimensions) which would fill in the last bit of the gap.  It
wouldn't be an isomorphism, but we have said that we can't have an
isomorphism so this is the next best thing.

Thus both methods would lead to a Dirac operator.  This raises two
questions: are the Dirac operators isomorphic?  and which is the best
method?  The answers are: ``yes'' and ``the cotangent method''.  The
first one comes from the fact that the spinor bundles are constructed
using \emph{completions} of the tangent and cotangent bundles with
respect to the inner product topology, rather than the bundles
themselves, and these completions are isomorphic.

Since the operators are equivalent, it may seem surprising that we
therefore claim that one method is superior to the other.  The
reasoning is simple: the tangent method uses the inner products on
both the tangent and cotangent bundles (the former in the
construction, the latter in filling the gap) whereas the cotangent
method only uses the inner product on the cotangent bundle (which is
no longer induced by that on the tangent bundle).  The two inner
products are now independent and therefore given a choice between
using both or using only one, we lean towards the simpler option.

Before proceeding, we note that one option when encountering problems
of this nature in loop spaces is to alter the type of loop used.
Certainly, using something like \(H^1\)-Sobolev loops would define a
Hilbert manifold of loops and thus the inner product on the tangent
bundle would identify the tangent and cotangent bundles as in finite
dimensions.  However, the remarkable isomorphism would then fail and
so we would be looking for a way to construct a map which on fibres
looks like: \(\m{L}(H_1, H_2) \to H_1^* \wotimes H_2\).  It
may not seem so, but this is exactly the same type of problem as we
have above: the space \(\m{L}(H_1, H_2)\) is isomorphic to the dual of
\(H_1^* \wotimes H_2\) and so we are looking for a map from
a space to its dual when we already have a map the other way around.
Therefore we gain nothing by altering the type of loop.

\subsection{Inner Products and Hilbert Completions}
\label{sec:introhilbert}

In this section we pick up on a remark made above.  In the previous
section it was stated that the two constructions of the Dirac operator
in infinite dimensions produce equivalent operators because the
construction depends on the completions of the bundles with respect to
the inner product topology, rather than the bundles themselves.  It is
this property that gives a little more substance to the study of inner
products on infinite dimensional vector bundles.

The question of \emph{existence} of an inner product on an infinite
dimensional vector bundle is solved in a similar manner to that in
finite dimensions.  We need two conditions to be satisfied, one on the
base space and one on the typical fibre:
\begin{enumerate}
\item The base manifold is \emph{smoothly paracompact}, in that it
  admits smooth partitions of unity.

\item The typical fibre admits inner products%
\footnote{%
  This is not a trivial condition.  The direct product of a countable
  number of copies of \R does not admit any inner products.
}%
.
\end{enumerate}
Providing these two are satisfied we can define an inner product
exactly as in finite dimensions by picking local inner products and
summing them using a partition of unity.

Thus mere existence is not a problem and one might feel that
\pageref{page:last}~pages is a little long for a discussion as to why
one particular inner product is better than any other.  The truth of
the matter is that one wants more than just an inner product.  What is
needed is that the fibrewise completions of the cotangent bundle with
respect to the inner product fit together to yield a bundle of Hilbert
spaces.  If one starts with an arbitrary inner product there is no
guarantee that this will happen.

To make this specific, we recall the definition of equivalent inner
products:
\begin{defn}
  Let \(V\) be a locally convex topological vector space.  Let
  \(\ipc_1\) and \(\ipc_2\) be two inner products on \(V\) with
  corresponding norms \(\norm_1\) and \(\norm_2\).  We say that these
  inner products are \emph{equivalent} if the norms are equivalent.
  That is, there are constants \(a, b > 0\) such that \(a \norm[v]_1
  \le \norm[v]_2 \le b \norm[v]_1\) for all \(v \in V\).
\end{defn}

The following results are standard from Banach space theory:
\begin{lemma}
  Let \(H_1\) and \(H_2\) be the completions of \(V\) with respect to
  \(\norm_1\) and \(\norm_2\) respectively.  Then \(\ipc_1\) and
  \(\ipc_2\) are equivalent if and only if the identity map on \(V\)
  extends to an isomorphism \(H_1 \cong H_2\).

  Let \(\ipc\) be an inner product on \(V\) with Hilbert completion
  \(H\).  Let \(g \in \gl(V)\).  Define \(\ipc_g\) by \(\ip{u}{v}_g =
  \ip{g u}{g v}\).  Then \(\ipc\) is equivalent to \(\ipc_g\) if and
  only if \(g\) extends to an operator in \(\gl(H)\).
\end{lemma}

From this it is clear that for the fibrewise Hilbert completions to
form a bundle then the equivalence class of the inner product must be
constant.  This leads to four types of inner product which we define
in terms of an associated principal bundle.  We think of a point in
this principal bundle as being an isomorphism from the corresponding
fibre to the model space.  In infinite dimensions it is rare to use
the full general linear group as this is either not a Lie group or is
contractible.
\begin{defn}
  We classify the inner products on a vector bundle according to how
  many of the following statements are satisfied.

  \begin{enumerate}
  \item The \emph{basic} inner product: the vector bundle admits a
    smooth choice of inner product on its fibres.

  \item The \emph{completable} inner product: the vector bundle admits
    a smooth choice of inner product on its fibres which map to a
    fixed equivalence class under the action of the principal bundle.

  \item The \emph{weakly locally trivial} inner product: the vector
    bundle admits a completable inner product and the principal bundle
    can be altered by a homotopy so that the inner product is mapped
    to a fixed inner product under its action.

  \item The \emph{locally trivial} inner product: the vector bundle
    admits a completable inner product which maps to a fixed inner
    product under the action of the principal bundle.
  \end{enumerate}
\end{defn}

Now that we have explicitly introduced a principal bundle, we can
rephrase these definitions in terms of the action of the group on the
vector space.  This will make clearer what is meant by a ``weakly
locally trivial'' inner product.  Using the fact that an inner product
on a space is the same thing as an inclusion with dense image into a
Hilbert space with an inner product%
\footnote{%
  Unless otherwise stated, when equipping a Hilbert space with an
  inner product we shall assume that it generates the given topology.
}%
, we can also phrase the statements using Hilbert spaces rather than
inner products.

\begin{proposition}
  Let \(G\) be a Lie group acting on a vector space \(V\).  Let \(P
  \to X\) be a principal \(G\)-bundle over a manifold \(X\) (which we
  assume to be smoothly paracompact).  Let \(E := P \times_G V\) be
  the associated vector bundle.   The following four conditions are
  equivalent, in order, to the different types of inner product given
  above:
  \begin{enumerate}
  \item \(V\) admits an inner product; equivalently, there is an
    inclusion \(V \to H\) with dense image of \(V\) into a Hilbert
    space, \(H\).

  \item The action of \(G\) on \(V\) preserves an equivalence class of
    an inner product on \(V\); equivalently, \(G\) acts on the diagram
    \(V \to H\) but not necessarily by isometries.

  \item The action of \(G\) on \(V\) preserves an equivalence class of
    an inner product on \(V\) and there is a subgroup \(K\) of \(G\)
    homotopic to \(G\) with an action on \(V\) by isometries such that
    this action is homotopic to the action which factors through
    \(G\); equivalently, the induced action of \(K\) on the diagram
    \(V \to H\) can be altered by homotopy so that it acts by
    isometries.

  \item The action of \(G\) on \(V\) is by isometries with respect to
    a fixed inner product; equivalently, the action of \(G\) on the
    diagram \(V \to H\) is by isometries.
  \end{enumerate}
\end{proposition}

We can now state the main theorem of this paper:
\begin{theorem}
  The cotangent bundle of the loop space of a smooth manifold
  considered as a bundle with structure group \(L \gl_n(\R)\) admits a
  weakly locally trivial inner product.

  The cotangent bundle of the loop space of a Riemannian manifold
  considered as a bundle with structure group \(L O_n\) does not admit
  a locally trivial inner product.
\end{theorem}

Compare this with the well-known analogous theorem for the tangent
bundle:
\begin{theorem}
  The tangent bundle of the loop space of a smooth manifold
  considered as a bundle with structure group \(L \gl_n(\R)\) admits a
  weakly locally trivial inner product.

  The tangent bundle of the loop space of a Riemannian manifold
  considered as a bundle with structure group \(L O_n\) admits a
  locally trivial inner product.
\end{theorem}

This inner product being given by the formula:
\[
\ip{\alpha}{\beta}_\gamma = \int_{S^1}
\ip{\alpha(t)}{\beta(t)}_{\gamma(t)} d t,
\]
where we identify the tangent space of the loop space with the loop
space of the tangent space.

The use of principal bundles has a significant advantage over just
writing down formul\ae\ such as the one above.  When writing down a
formula one then has to go to considerable lengths to prove any local
triviality statements that one may wish to use whereas the local
triviality follows naturally if everything is done using principal
bundles.

This seems an appropriate point to mention one aspect of the theory
that influences the flavour of the discussion without introducing any
change in the mathematics.  As part of our quest we shall have to pick
a Hilbert completion of \((L \R)^*\).  Because all the spaces involved
are reflexive, we can equally choose a dense Hilbert subspace of \(L
\R\).  Since it is conceptually easier to visualise subspaces of \(L
\R\) than superspaces of \((L \R)^*\), we tend to work in this dual
picture.  Reflexivity ensures that we introduce no complications by
doing so.

\subsection{The Inner Product on the Cotangent Bundle - an Overview}
\label{sec:introcotangent}

We shall now give an overview of the construction of the inner product
on the cotangent bundle.  The method that we employ is to look at the
structure group.

From the point of view of algebraic topology the construction is very
simple.  There is no equivalence class of inner product on the model
space of the cotangent bundle that is preserved by the action of \(L
\gl_n(\R)\), or even by \(L O_n\), but the polynomial loop group,
\(L_\pol O_n\) preserves many equivalence classes.  Since the
polynomial loop group is homotopic to the smooth loop group, see for
example~\cite{apgs}, we choose a reduction of the structure group from
the smooth loop group to the polynomial loop group.  There is a little
work to show that the action of the polynomial loop group can be
altered through homotopies to one by isometries with respect to some
fixed inner product, but this is not hard.  Thus we have a weakly
locally trivial inner product.

The purpose of the rest of the \pageref{page:last}~pages of this paper
is to reduce the number of choices in that last paragraph to a minimum
and to make them as global as possible.  Ultimately, we end up with
one global choice which is essentially the reference inner product on
the model space of the cotangent bundle.  Our input to the machinery
is a Riemannian manifold -- which is the same input needed to define
the inner product on the tangent bundle.  However, our method of
construction is somewhat more complicated.

We shall explain the method for the space of based loops, \(\Omega
M\).  This will enable us to get to the central idea without too many
details.  We assume that \(M\) is simply-connected so that its loop
spaces are connected.  This implies that \(M\) is orientable.

The Riemannian structure on \(M\) defines the Levi-Civita connection.
This in turn defines the holonomy operator, \(h : \Omega M \to S
O_n\).  Now \(S O_n\) is the classifying space of the group \(\Omega S
O_n\) and the holonomy operator is a classifying map for the principal
\(\Omega S O_n\)-bundle associated to the tangent bundle of \(\Omega
M\) -- and thus also to the cotangent bundle.  The parallel transport
operator defines an explicit isomorphism from the (co)tangent bundle
to the corresponding pull-back bundle.

Since the polynomial loop group, \(\Omega_\pol S O_n\), is homotopic 
to the smooth loop group, \(\Omega S O_n\), the classifying spaces are
the same.  Thus there is a principal \(\Omega_\pol S O_n\)-bundle over
\(S O_n\) which includes into the natural \(\Omega S O_n\)-bundle.
Therefore, using the holonomy and parallel transport maps, we get a
principal \(\Omega_\pol S O_n\)-bundle over \(\Omega M\) which is a
natural subbundle of the principal \(\Omega S O_n\)-bundle associated
to the tangent and cotangent bundles.

As stated above, the action of \(\Omega_\pol S O_n\) on the model
space of the cotangent preserves an equivalence class of an inner
product, though does not act by isometries.  The choice of this
equivalence class -- for there are many -- is part of our one choice.
We can therefore define the Hilbert completion of the cotangent
bundle.  The action of \(\Omega_\pol S O_n\) can be gently altered to
one of isometries whereupon we get an inner product.  The choice of
this inner product is the other part of our one choice (we count these
as one choice since the inner product determines its equivalence
class).

When considering the full loop space we use the fact that there is a
locally trivial fibration \(\Omega M \to L M \to M\) and essentially
repeat the above construction fibre-by-fibre on \(L M \to M\).

We give an explicit construction of the \(\Omega_\pol S O_n\)-bundle
over \(S O_n\) (actually, we construct an \(L_\pol S O_n\)-bundle) and
show that it is locally trivial.  Thus although the homotopy
equivalence \(\Omega_\pol S O_n \to \Omega S O_n\) is part of the
background of the construction, we never actually use it.  In fact,
the bundle that we construct shows that the homotopy groups of
\(\Omega S O_n\) are a direct summand of those of \(\Omega_\pol S
O_n\): the existence of this bundle means that there is a map \(B
\Omega S O_n \to B \Omega_\pol S O_n\).  Since the \(\Omega_\pol S
O_n\)-bundle includes naturally in the \(\Omega S O_n\)-bundle, the
composition of the above map on classifying spaces with the natural
map \(B \Omega_\pol S O_n \to B \Omega S O_n\) is homotopic to the
identity.  Ralph Cohen has suggested that further study of this bundle
might yield an alternative proof of the homotopy equivalence of
\(\Omega_\pol S O_n\) with \(\Omega S O_n\), however that is beyond
the scope of this paper.

\subsection{Acknowledgements and History}
\label{sec:introack}

The central idea of this paper -- the construction of the polynomial
loop bundle -- places this paper as the latest in a loosely defined
series: \cite{jm2}, \cite{rcas}, and \cite{as2}.  In the first of
these, Morava attempted to construct an isomorphism for an almost
complex manifold \(M\) between the tangent bundle of the loop space,
\(L T M\), and a bundle of the form \(e_1^* T M \otimes_\C L \C\).
Here, \(e_1 : L M \to M\) is the map which evaluates a loop at time
\(1\) and every bundle is considered to be complex.  The argument
broke down at one crucial step and the papers \cite{rcas} and
\cite{as2} grew out of considering the question as to when that
crucial step could be made to work.  This was found to be highly
restrictive and implied, for example, that the tangent bundle of the
based loop space of \(M\) was trivial.

One consequence which would follow from the existence of an
isomorphism \(L T M \cong e_1^* T M \otimes_\C L \C\) would be the
existence of a sub-bundle modelled on the polynomial loop space.  In
fact, for any class of loops there would be a bundle with the
appropriate fibre constructed as \(e_1^* T M \otimes_\C L^\alpha \C\).
Close examination of \cite{jm2} reveals that Morava's method was
essentially to construct the polynomial loop bundle fibrewise.  His
mistake was to assume that from this one could globally pick-out a
finite dimensional sub-bundle.

The point of view of this paper is that the polynomial loop bundle is
as far as one needs to go -- as well as being as far as one can go.
We proffer the construction of the Dirac operator as evidence for
this.

The author would like to thank Rafe Mazzeo, Ralph Cohen, and Eldar
Straume for helpful conversations and to acknowledge the encouragement
of Jack Morava.

\subsection{Structure of the Paper}
\label{sec:intropaper}

This paper is structured as follows: 
\begin{description}
\item[Section~\nmref{sec:notation}:] In this section we gather in one
  place all the non-standard or unusual notation that we shall use in
  this paper.  This section is more than a reference section in that
  notation defined here will not be formally defined elsewhere.

\item[Section~\nmref{sec:polgrp}:] In this section we study the
  polynomial loop groups and construct the universal polynomial loop
  bundles over the classifying spaces.  For technical reasons -- which
  are given -- we concentrate on the cases \(U_n\), \(S U_n\), and \(S
  O_n\).

\item[Section~\nmref{sec:dual}:] In this section we use the work of
  section~\ref{sec:polgrp} to construct an inner product on the dual
  of a loop bundle and the associated Hilbert bundle.  We also
  consider the properties of this construction and the relation of
  loop bundles to twisted K-theory.

\item[Section~\nmref{sec:diracloop}:] In this section we construct a
  Dirac operator over the loop space of a string manifold using the
  inner product on the cotangent bundle.  This section also contains a
  brief summary of the main results on infinite dimensional spin.

\item[Appendix:] This contains some results that may be of interest
  about the general problem of inner products on the space of
  distributions.
\end{description}

\newpage
\section{Notation}
\label{sec:notation}

This paper is somewhat heavy on notation.  Therefore, we have included
this section here as a reference point for the bemused reader.  Here
we have collected together the notation for all the reasonably
standard objects that we use.  The following definitions have not been
included here as they are the main subject of study in various
sections of this paper:
\begin{enumerate}
\item The polynomial loop groups: \(\Omega_\pol G\) and \(L_\pol G\).
  Section~\ref{sec:polloops}.

\item The periodic and polynomial path spaces: \(P_\per G\) and
  \(P_\pol G\).  Section~\ref{sec:path}.

\item The periodic and polynomial vector bundles: \(P_\per V\) and
  \(P_\pol V\).  Section~\ref{sec:polvect}.

\item The polynomial loop bundles: \(L_\pol E\), \(L_\pol Q\), and
  \(L_\pol Q^\ad\).  Section~\ref{sec:poly}.

  Also defined in section~\ref{sec:poly} are various spaces used in
  the construction of the polynomial loop bundles.  As these are not
  used elsewhere we shall not list them here.
\end{enumerate}

We have tried to choose notation that is as clear as possible by
choosing notation that is relatively bracket free.  The issue is
further complicated by the fact that this topic mixes geometry and
functional analysis.  Notation that is clear when geometrically viewed
may not be so from the point of view of a functional analysts.  As the
intended audience consists primarily of geometers, we have gone for
clarity in the geometrical viewpoint, with apologies to any functional
analysts that may be present.

\subsection{The Circle}
\label{sec:notecirc}

In this paper we have two views of the circle.  One is as the domain
of loops, the other as a Lie group.  We regard loops as periodic paths
from \R and thus wish to identify the domain of loops with \(\R/\Z\).
When thinking of the circle as a Lie group, we think of it as \(U_1\)
sitting inside \(M_1(\C) = \C\).  We shall use the notation \(S^1\)
for \(\R/\Z\) and \T for \(U_1\).  We shall write \(t\) for the
parameter in \(S^1\) and \(z\) in \T, with relationship \(z = e^{2 \pi
i t}\).

\subsection{Loop and Path Spaces of Fibre Bundles}
\label{sec:notepath}

Let \(M\) be a finite dimensional smooth manifold.  We shall write \(L
M\) for the manifold of smooth maps \(S^1 \to M\) and \(P M\) for the
manifold of smooth maps \(\R \to M\).  Since we are viewing a loop as
a periodic path with period \(1\), \(L M\) is a submanifold of \(P
M\).  By regarding \(M\) as the space of constant paths, we view \(M\)
as a submanifold of \(L M\), whence also of \(P M\).

Let \(F \to X \to M\) be a locally trivial fibre bundle which is
either a vector bundle, a principal bundle, or a bundle of Lie groups.
The loop space of \(X\) is a locally trivial fibre bundle over \(L
M\).  If \(X\) is \emph{orientable} -- that is, trivialisable over any
loop -- then the fibre is \(L F\); otherwise it will vary on the
components of \(L M\).  In the situations encountered in this paper
the bundles will always be orientable.  In all cases, \(P X
\to P M\) is a locally trivial fibre bundle with fibre \(P F\).

We shall define various pull-backs of the bundles \(L X \to L M\) and
\(P X \to P M\).  The guide to our notation is that we shall label the
pull-backs by adjoining appropriate superscripts to the \(L\) or
\(P\).  The convention will be to read from left to right: that is,
the leftmost label happened first.  Thus \(L^{a, b} X\) denotes the
bundle \(L X\) pulled back via first \(a\) and then \(b\).

We shall label the fibre of a bundle over a particular point by
adjoining the label of point as a subscript to the \(L\) or \(P\).
When the \(L\) or \(P\) is decorated by an additional subscript, say
\ttype, the fibre label will be to the right of this.  Thus
\(L_{\type, \gamma} X\) is the fibre of \(L_\type X\) over \(\gamma
\in L M\) (the additional subscripts are \tpol and \tper which will be
defined in section~\ref{sec:polgrp}).

We shall now describe the various pull-back bundles that we shall use:
\begin{enumerate}
  \item \(P^L X\) is the pull-back (or restriction) of \(P X\) to \(L
    M\); thus for \(\gamma \in L M\), \(P^L_\gamma X = P_\gamma X\).
    Note that \(L X\) is a sub-bundle of \(P^L X\).

  \item \(P^M X\) and \(L^M X\) are the pull-backs of, respectively,
    \(P X\) and \(L X\) to \(M\).  Again, \(P^M_p X = P_p X\).

    Note that as \(p\) is here regarded as a constant path, the paths
    (resp.\ loops) in \(P_p X\) (resp.\ \(L_p X\)) lie above a single
    point in \(M\).  Thus they lie in a single fibre of \(X\).  Hence
    \(P_p X = P(X_p)\) (resp.\ \(L_p X = L(X_p)\)).

  \item For \(t \in \R\), let \(e_t : P M \to M\) be the map which
    evaluates a path at time \(t\).  Let \(X^t \to P M\) be the
    pull-back of \(X\) via \(e_t\).  We shall use the same notation
    for \(X^t\) restricted to \(L M\) to avoid too many superscripts.
    Likewise, let \(P^{M,t} X\) and \(L^{M,t} X\) be the pull-backs of
    \(P^M X\) and \(L^M X\), respectively, via \(e_t\).  Thus
    \(X^t_\gamma = X_{\gamma(t)}\), \(P^{M,t}_\gamma X =
    P(X_{\gamma(t)})\), and \(L_\gamma^{M,t} X = L( X_{\gamma(t)})\).

    For each \(t \in \R\), there are evaluation maps \(P X \to X^t\)
    which we denote again by \(e_t\).  Over \(L M\), we have the
    identity: \(X^{t + 1} = X^t\).
\end{enumerate}

\subsection{Function Spaces and Function Bundles}
\label{sec:notefun}

The function spaces that we shall use in this paper will be spaces of
maps from \(S^1\) to some finite dimensional real or complex vector
space.  We shall base our notation on that from differential geometry
rather than functional analysis and use a similar convention to that
above.  Thus a space of maps from \(S^1\) will be denoted by an \(L\)
decorated in some fashion.

Since we are using \(L\) to denote maps from the circle into some
vector space, we shall use \(\m{L}(X,Y)\) for continuous linear maps
from one topological vector space to another.  Where the target space
is the same as the source, we shall abbreviate this to \(\m{L}(X)\).

We shall now describe the various function spaces that we shall use in
terms of maps from the circle into \C.  For maps into \(\C^n\), we
tensor with \(\C^n\): thus \(L \C^n = L \C \otimes \C^n\); for maps
into \(\R^n\), we take the underlying real space of the maps into
\(\C^n\).

\begin{enumerate}
\item \(L \C\): smooth maps.

\item \(L^2 \C\): square-integrable maps.

\item \(L_\pol \C = \C[z^{-1}, z]\): Laurent polynomials in \C.

\item \(L^* \C\): distributions -- the dual of \(L \C\).

\item \(L^{2,*} \C\): the dual of \(L^2 \C\).

\item \(L_r \C\), \(r > 1\): smooth maps which
  extend holomorphically over an annulus of outer radius \(r\) and
  inner radius \(r^{-1}\) and are smooth on the boundary.

\item \(L_r^2 \C\), \(r > 1\): smooth maps which extend
  holomorphically over an annulus of outer radius \(r\) and inner
  radius \(r^{-1}\) and are square-integrable on the boundary.

\item \(L_r^* \C\), \(r > 1\): smooth maps which extend
  holomorphically over an annulus of outer radius \(r\) and inner
  radius \(r^{-1}\) and are distributions on the boundary.

\item \(L_r^{2,*} \C\), \(r > 1\): the dual of \(L_r^2 \C\).

\item \(L_r^{(2,*)} \C\), \(r > 1\): smooth maps which extend
  holomorphically over an annulus of outer radius \(r\) and inner
  radius \(r^{-1}\) and are dual to square-integrable on the boundary.
\end{enumerate}

The last space is, of course, just \(L_r^2 \C\).  However, we are
viewing it as the image in \(L^{2,*} \C\) of \(L_r^2 \C\) under the
conjugate linear isomorphism \(L^2 \C \to L^{2,*} \C\).

The penultimate space in the above list has an interesting
interpretation.  Within the space of formal power series,
\(\C[\![z^{-1}, z]\!]\), one can consider those power series that
converge on a \emph{formal} annulus of \emph{outer} radius \(r^{-1}\)
and \emph{inner} radius \(r\) for some \(r > 1\) and satisfy some
condition on the boundary.  It is not hard to show that this space is
(conjugate) dual to some space of the form \(L_r^a \C\) for some
appropriate boundary condition.  Thus \(L_r^{2,*} \C\) is conjugate
dual to \(L_{r^{-1}}^2 \C\).  One consequence of this interpretation
is the following identity:
\[
(L_r^{2,*})_r \C = (L_{r^{-1}}^2)_r \conj{\C} = L^2 \conj{\C} =
L^{2,*} \C.
\]
The crucial step here is the observation that the annuli of radii
\((r,r^{-1})\) and of radii \((r^{-1}, r)\) cancel out.

Let \(E \to M\) be a finite dimensional orientable vector bundle over
a finite dimensional smooth manifold.  The loop space, \(L E\), is an
infinite dimensional vector bundle over \(L M\) modelled on \(L \F^n\)
for some \(n\), where \F is either \R or \C.  We shall consider
various related bundles where we modify the fibre of \(L E\) from \(L
\F^n\) to some other function space.  We shall label these by
decorating the \(L\) as for the function spaces above.

The standard ones are \(L^2 E\) and \(L^* E\) having fibre \(L^2
\F^n\) and \(L^* \F^n\) respectively.  The bundle \(L^* E\) is the
dual of \(L E\): a fibre, \(L^*_\gamma E\), is the space of continuous
linear maps \(L_\gamma E \to \F\).  The simplest way to define \(L^2
E\) is to observe that the action on \(L \F^n\) of the structure group
of \(L E\), namely \(L \gl_n(\F)\), extends to an action on \(L^2
\F^n\).  Thus \(L^2 E\) is constructed from the principal bundle of
\(L E\) in the usual way.  Fibrewise, it can be viewed as the Hilbert
completion of \(L E\) with respect to the inner product:
\[
\ip{\alpha}{\beta}_\gamma := \int_{S^1}
\ipv{\alpha(t)}{\beta(t)}_{\gamma(t)} d t
\]
where \(\ipvc\) is some smooth choice of inner product on the fibres
of \(E\).  With this approach, one needs to show that this fibrewise
completion does result in a locally trivial Hilbert bundle.

The definitions of the bundles \(L^2_r E\), \(L^*_r E\), and the
others is the core of this paper.  They will turn out to be locally
trivial bundles modelled on the corresponding function spaces.

Finally, it is a standard fact from the differential topology of loop
spaces that \(T L M = L T M\) but that \(T^* L M \ne L T^* M\).  For
the second, observe that in the case of \(\R^n\), \(T^* L \R^n = L^*
\R^n \times L \R^n\) but \(L T^* \R^n = L \R^n \times L \R^n\).  Thus
\(T\) and \(L\) commute whilst \(T^*\) and \(L\) do not.  The notation
we have introduced above provides another way of writing the cotangent
bundle, namely \(L^* T M\).  With this notation, \(T\) and \(L\)
continue to behave well since \(T^* L M = L^* T M\).

To continue into absurdity, note that \(L^*\) and \(T\) do not commute
even when \(L^* M\) makes sense (i.e.~when \(M\) is \(\R^n\)) since
\(T L^* \R^n = L^* \R^n \times L^* \R^n\) and \(L^* T \R^n = L^* \R^n
\times L \R^n\).  To break the bounds of absurdity and enter in to the
ridiculous, observe that \(T^* L^* \R^n = L \R^n \times L^* \R^n = T^*
L \R^n = L^* T \R^n\).  Thus our identities are: \(T L = L T\), \(T^*
L = L^* T = T^* L^*\), and \(T L^* = (T L)^*\).

\newpage
\section{The Polynomial Loop Group}
\label{sec:polgrp}

In this section we consider the group of polynomial loops in a
compact, connected Lie group.  This was studied extensively
in~\cite{apgs} with some further work appearing in~\cite{gs2} in the
case of \(U_n\).  We start with some general results on polynomial
loops before constructing the \(L_\pol G\)-bundle over \(G\) for \(G\)
each of \(U_n\), \(S U_n\), and \(S O_n\).  We conclude by
constructing the corresponding vector bundles.

\subsection{Polynomial Loops}
\label{sec:polloops}

The definition of the polynomial loop group appears in~\cite[\S
3.5]{apgs}.  We repeat that definition here.

\begin{defn}
  \label{def:polgrp}
  Let \(G\) be a compact, connected Lie group.  Fix an embedding of
  \(G\) as a subgroup of \(U_n\) for some \(n\).  This exhibits \(G\)
  as a submanifold of \(M_n(\C)\).  The \emph{polynomial loop group}
  of \(G\), \(L_\pol G\), is defined as the space of those loops in
  \(G\) which when expanded as a Fourier series in \(M_n(\C)\) are
  finite Laurent polynomials.  The group of based loops, \(\Omega_\pol
  G\), is the subgroup of \(L_\pol G\) of loops \(\gamma\) with
  \(\gamma(0) = 1_G\).
\end{defn}

\begin{remark}
  The following comments appear in~\cite[\S 3.5]{apgs}:

  \begin{enumerate}
  \item The choice of the embedding of \(G\) in \(U_n\) is immaterial.

  \item The space \(L_\pol G\) is the union of the subspaces
    \(L_{\pol,N} G\) consisting of those loops with Fourier series of
    the form:
    \[
    \sum_{k = -N}^N \gamma_k z^k.
    \]
    These spaces are naturally compact.  The topology on \(L_\pol G\)
    is the direct limit topology of this union.

  \item The free polynomial loop group is the semi-direct product of
    the based polynomial loop group and the constant loops.

  \item The group \(L_\pol G\) does not have an associated Lie
    algebra, although the Lie algebra \(L_\pol \mf{g}\) is often
    linked to it.

  \item If \(G\) is semi-simple then \(L_\pol G\) is dense in \(L
    G\).

  \item In the case of the circle, \(\Omega_\pol S^1 = \Z\) and so
    \(L_\pol S^1 = S^1 \times \Z\).
  \end{enumerate}
\end{remark}

The following is \cite[proposition~8.6.6]{apgs}:

\begin{proposition}
  The inclusion \(\Omega_\pol G \to \Omega G\) is a homotopy
  equivalence.
\end{proposition}

Since \(L_\pol G \cong \Omega_\pol G \times G\) and \(L G \cong \Omega
G \times G\) as spaces (although not generally as groups), this holds
for the unbased loops as well.

Although the definition of \(L_\pol G\) does not depend on the
embedding of \(G\) in \(U_n\), it is useful to have such an embedding
to investigate the structure of \(L_\pol G\) in a little more detail.
We consider loops of the form \(t \to \exp(t \xi)\) for suitable \(\xi
\in \mf{g}\).  The main result is the following:

\begin{proposition}
  \label{prop:liepol}
  Let \(G\) be a compact, connected Lie group, \(\mf{g}\) its Lie
  algebra.  For \(\xi \in \mf{g}\), let \(\eta_\xi : \R \to G\) denote
  the path \(\eta_\xi (t) = \exp(t \xi)\).

  Let \(\xi_1, \xi_2 \in \mf{g}\) be such that \(\exp(\xi_1) =
  \exp(\xi_2)\).  Then \(\eta_{-\xi_1} \eta_{\xi_2}\) is a polynomial
  loop in \(G\).
\end{proposition}

As part of the proof of this, we shall prove the following useful
result for the unitary group:

\begin{lemma}
  \label{lem:comlie}
  Let \(g \in U_n\).  There exists \(\zeta \in \exp^{-1}(g) \subseteq
  \mf{u}_n\) such that \([\zeta, \xi] = 0\) for all \(\xi \in
  \exp^{-1}(g)\).
\end{lemma}

The proofs of these rely on the simple structure in \(U_n\) of the
centraliser of any particular element.  For \(g \in U_n\), define
\(C(g)\) and \(Z(g)\) to be the centraliser of \(g\) and its centre.
That is, \(C(g) := \{h \in G : h^{-1} g h = g\}\) and \(Z(g) =
Z(C(g))\).  Clearly, \(g \in Z(g)\).

\begin{lemma}
  For any \(g \in U_n\), \(Z(g)\) is a torus.
\end{lemma}

\begin{proof}
  The group \(C(g)\) is a closed subgroup of \(U_n\), hence its centre
  is a closed abelian subgroup of \(U_n\).  In particular, it is
  compact.  Therefore, it is a torus if and only if it is connected.

  Recall that two diagonalisable matrices commute if and only if they
  are simultaneously diagonalisable.  This condition does not rely on
  the eigenvalues of either matrix but only on the eigenspaces.

  Let \(h \in Z(g)\).  As \(h\) is unitary, it is orthogonally
  diagonalisable.  Let \(\lambda_1, \dotsc, \lambda_l\) be the
  distinct eigenvalues of \(h\) with associated eigenspaces \(E_1,
  \dotsc, E_l\).  For each \(j\), let \(s_j \in [- i \pi, i \pi) \) be
  such that \(e^{s_j} = \lambda_j\).

  Define \(\alpha : [0, 1] \to U_n\) to be the path such that
  \(\alpha(t)\) has eigenvalues \(e^{t s_j}\) and corresponding
  eigenspaces \(E_j\).  Then \(\alpha(0) = 1_n\) and \(\alpha(1) = h\)
  so \(\alpha\) is a path from \(1_n\) to \(h\).  By construction,
  \(\alpha(t)\) for \(t \ne 0\) has the same eigenspaces as \(h\) and
  therefore \(\alpha(t)\) commutes with exactly the same elements of
  \(U_n\) that \(h\) commutes with.  Hence as \(h \in Z(g)\),
  \(\alpha(t) \in Z(g)\).
\end{proof}

\begin{proof}[Proof of lemma~\nmref{lem:comlie}]
  As \(Z(g)\) is a torus, it is a connected compact Lie group.
  Therefore, the exponential map is surjective and so there is some
  \(\zeta \in \mf{z}(g) \subseteq \mf{u}_n\) with \(\exp(\zeta) =
  g\).  As \(\zeta \in \mf{z}(g)\), \(\exp(t \zeta) \in Z(g)\) for all
  \(t \in \R\).

  Let \(\xi \in \mf{u}_n\) be such that \(\exp(\xi) = g\).  Then for
  all \(t \in \R\), \(\exp(t \xi)\) commutes with \(g\).  Hence
  \(\exp(t \xi) \in C(g)\) for all \(t\).  Thus \(\exp(t \zeta)\) and
  \(\exp(t' \xi)\) commute for all \(t, t' \in \R\).  Hence \([\zeta,
  \xi] = 0\).
\end{proof}

Using this we can prove proposition~\ref{prop:liepol}.

\begin{proof}[Proof of proposition~\nmref{prop:liepol}]
  Firstly, note that it is sufficient to prove this in the case of the
  unitary group.  For if \(\eta_{-\xi_1} \eta_{\xi_2}\) is a loop in
  \(G\) which is a polynomial loop when \(G\) is considered as a
  subgroup of \(U_n\) then, by definition, \(\eta_{-\xi_1}
  \eta_{\xi_2}\) is a polynomial loop in \(G\).

  Secondly, note that it is sufficient to consider the case where
  \(\xi_2 = 0\).  This forces \(\exp(\xi_1) = 1_n\).  To deduce the
  general case from this simpler one, note that by
  lemma~\ref{lem:comlie} that there is some \(\zeta \in \mf{u}_n\)
  with \(\exp(\zeta) = \exp(\xi_1)\) (whence also \(\exp(\xi_2)\))
  such that \([\zeta,\xi_j] = 0\).  Then \(\exp(\xi_j - \zeta) = 1_n\)
  so, by assumption, \(\eta_{(\xi_j - \zeta)}\) is a polynomial loop.
  The identity:
  \[
  \eta_{-\xi_1} \eta_{\xi_2} = \eta_{-\xi_1} \eta_\zeta \eta_{-\zeta}
  \eta_{\xi_2} = \eta_{(-\xi_1 + \zeta)} \eta_{(\zeta - \xi_2)}.
  \]
  demonstrates that this is a polynomial loop.

  Thus we need to show that \(\eta_\xi\) is a polynomial loop if
  \(\exp(\xi) = 1\).  To show this, we diagonalise \(\xi\).  If \(s\)
  is an eigenvalue of \(\xi\) then \(e^{s}\) is an eigenvalue of
  \(\exp(\xi) = 1\).  The eigenvalues of \(\xi\) therefore lie in \(2
  \pi i \Z\).  Hence there is a basis of \(\C^n\) with respect to
  which \(\eta_\xi\) is the path:
  \[
  t \to
  \begin{bmatrix}
    e^{2 \pi i t k_1} & 0 & \dots & 0 \\
    0 & e^{2 \pi i t k_2} & \dots & 0 \\
    \hdotsfor{4} \\
    0 & 0 & \dots & e^{2 \pi i t k_n}
  \end{bmatrix}
  \]
  for some \(k_j \in \Z\).  Since \(e^{2 \pi i t k} = z^k\) for \(k
  \in \Z\), this is a polynomial loop (viewed as a periodic path).
\end{proof}

Note that for a general group \(G\), although the loop \(\eta_{\xi_1}
\eta_{-\xi_2}\) lies in \(\Omega_\pol G\), there may be no
factorisation in \(G\) as \(\eta_{\xi_1 - \zeta} \eta_{\zeta -
\xi_2}\) since in a general Lie group there may not be any \(\zeta \in
\mf{g}\) satisfying the required properties.

\subsection{The Path Spaces}
\label{sec:path}

In the light of the homotopy equivalence, \(\Omega_\pol G \simeq
\Omega G \simeq \Omega_{\text{cts}} G\), the classifying space of
\(\Omega_\pol G\) is (homotopy equivalent to) \(G\) itself.  Since
\(\Omega_\pol G\) acts on \(L_\pol G\), there is a natural \(L_\pol
G\)-principal bundle over \(G\).  In this section we shall give an
explicit construction of this bundle.  We shall also construct a
similar bundle for the smooth loop group.  These bundles will be
denoted by \(P_\pol G\) and \(P_\per G\) (the ``\(\per\)'' stands for
``periodic'').

To demonstrate that these are principal bundles with the appropriate
fibre we have to show two things: firstly, that the bundles are
locally trivial; and secondly, that the fibres have an action of the
appropriate loop group which identifies the fibre with that group.
The second of these is straightforward, the first is simple for the
smooth case but is surprisingly difficult for the polynomial loop
group.  We shall only consider the cases of \(U_n\), \(S U_n\), and
\(S O_n\).

\begin{defn}
  \label{def:perpol}
  Let \(G\) be a compact, connected Lie group, \(\mf{g}\) its Lie
  algebra.  We define \(P_\per G\) and \(P_\pol G\) as
  follows:
  \begin{enumerate}
  \item \(P_\per G\) is the space of smooth paths \(\alpha : \R \to
    G\) with the property that \(\alpha(t + 1)\alpha(t)^{-1}\) is
    constant.

  \item \(P_\pol G \subseteq P_\per G\) consists of those paths of the
    form \(\eta_\xi \gamma\) for some \(\xi \in \mf{g}\) and \(\gamma
    \in L_\pol G\).
  \end{enumerate}

  The projection map \(P_\per G \to G\) is given by \(\alpha \to
  \alpha(1) \alpha(0)^{-1}\).  Notice that when restricted to \(P_\pol
  G\), this maps \(\eta_\xi \gamma\) to \(\exp(\xi)\).
\end{defn}

Recall from section~\ref{sec:polloops} that for \(\xi \in \mf{g}\) the
path \(\eta_\xi : \R \to G\) is defined as the path \(t \to \exp( t
\xi)\).

Observe that a path in \(P_\per G\) is completely determined by its
values on the interval \([0,1]\).  The motivation for the given
definition of \(P_\per G\) (and of \(P_\pol G\)) is that of holonomy.

It will sometimes be useful to consider an element of \(P_\per G\) to
be a pair \((g, \alpha) \in G \times P G\) such that \(\alpha(t + 1) =
g\alpha(t)\).  Here \(P G\) is \emph{all} smooth paths \(\R \to G\).
Although \(g\) is completely determined by \(\alpha\), this viewpoint
makes it more explicit.

We shall now investigate the desired properties of these spaces.
Neither is a group (unlike the analogous continuous situation), but
the group \(G\) acts in two ways:

\begin{lemma}
  \label{lem:pathconj}
  The group \(G\) acts on \(P_\per G\) by two actions:
  \[
  g \cdot_m \alpha = g \alpha, \qquad g \cdot_c \alpha = g \alpha
  g^{-1}.
  \]
  These actions restrict to actions on \(P_\pol G\).  For both
  actions, the action of \(G\) on itself by conjugation makes the
  projection \(P_\per G \to G\) \(G\)-equivariant (hence also for
  \(P_\pol G \to G\)).
\end{lemma}

\begin{proof}
  Let \(g \in G\) and \(\alpha \in P_\per G\).  Both \(g \alpha\) and
  \(g \alpha g^{-1}\) are smooth paths in \(G\) so we only need to
  check the periodicity condition.  Let \(h = \alpha(t + 1)
  \alpha(t)^{-1}\).  Then:
  \begin{align*}
    (g \alpha)(t + 1) (g \alpha)(t)^{-1} &= g \alpha(t + 1)
    \alpha(t)^{-1} g^{-1} = g h g^{-1}. \\
    (g \alpha g^{-1}) (t + 1) (g \alpha g^{-1})(t)^{-1} &= g \alpha(t
    + 1) g^{-1} g \alpha(t)^{-1} g^{-1} \\
    &= g \alpha(t + 1) \alpha(t)^{-1} g^{-1} = g h g^{-1}.
  \end{align*}
  This also proves the statement about the induced action on \(G\).

  If \(\alpha \in P_\pol G\) then \(\alpha\) is of the form \(\eta_\xi
  \gamma\) for some \(\xi \in \mf{g}\) and \(\gamma \in L_\pol G\).
  Let \(h\) be either \(g^{-1}\) or \(1_G\).  Then \(g \eta_\xi \gamma
  h = \eta_{(\Ad_g \xi)} g \gamma h\).  As \(L_\pol G\) is closed under
  left and right multiplication by \(G\), this lies in \(P_\pol G\) as
  required.
\end{proof}

\begin{proposition}
  Define an action of \(L G\) on \(P_\per G\) by sending \((\alpha,
  \gamma) \in P_\per G \times L G\) to the path \(t \to \alpha(t)
  \gamma(t)\).  This action is well-defined and identifies the fibres
  of \(P_\per G \to G\) with \(L G\).  It restricts to an action of
  \(L_\pol G\) on \(P_\pol G\) and identifies the fibres of \(P_\pol G
  \to G\) with \(L_\pol G\).
\end{proposition}

\begin{proof}
  The path \(t \to \alpha(t) \gamma(t)\) is a smooth path \(\R \to G\)
  (considering \(\gamma\) as a periodic path).  We need merely check
  the periodicity condition.  Since \(\gamma(t + 1) = \gamma(t)\) for
  all \(t \in \R\), we have:
  \begin{align*}
    (\alpha \gamma) (t + 1) (\alpha \gamma)(t)^{-1} &= \alpha(t + 1)
    \gamma (t + 1) \gamma(t)^{-1} \alpha(t)^{-1} \\
    &= \alpha(t + 1) \alpha(t)^{-1}.
  \end{align*}
  Hence \(\alpha \gamma \in P_\per G\).  This also shows that \(\alpha
  \gamma\) lies in the same fibre as \(\alpha\).

  For an inverse, let \(\alpha, \beta \in P_\per G\) be such that
  \(\alpha(1) \alpha(0)^{-1} = \beta(1) \beta(0)^{-1}\).  As
  \(\alpha\) and \(\beta\) lie in \(P_\per G\), this means that
  \(\alpha(t + 1)\alpha(t)^{-1} = \beta(t + 1) \beta(t)^{-1}\) for all
  \(t \in \R\).  Rearranging this yields \(\alpha(t + 1)^{-1} \beta(t
  + 1) = \alpha(t)^{-1} \beta(t)\).  Thus the path \(\gamma\) given by
  \(\gamma(t) = \alpha(t)^{-1} \beta(t)\) is a loop.  Moreover, it is
  smooth.  Clearly \(\alpha \gamma = \beta\) so this is the inverse
  map which identifies a non-empty fibre of \(P_\per G \to G\) with
  \(L G\).

  In the polynomial case, if \(\alpha \in P_\pol G\) and \(\gamma \in
  L_\pol G\) then by definition, \(\alpha = \eta_\xi \beta\)
  for some polynomial loop \(\beta\).  Therefore \(\alpha \gamma =
  \eta_\xi (\beta \gamma)\) and hence lies in \(P_\pol G\).
  
  Conversely, suppose that \(\alpha, \beta \in P_\pol G\) lie in the
  same fibre.  We need to show that the loop \(t \to \alpha^{-1}(t)
  \beta(t)\) is a polynomial loop.  Let \(\alpha = \eta_{\xi_1}
  \widehat{\alpha}\) and \(\beta = \eta_{\xi_2} \widehat{\beta}\)
  where \(\widehat{\alpha}\) and \(\widehat{\beta}\) are polynomial
  loops.  Since \(\alpha\) and \(\beta\) lie in the same fibre,
  \(\exp(\xi_1) = \exp(\xi_2)\).  Thus:
  \[
  \gamma = \widehat{\alpha}^{-1} \eta_{-\xi_1} \eta_{\xi_2}
  \widehat{\beta}.
  \]
  By proposition~\ref{prop:liepol}, the two terms in the centre give a
  polynomial loop, hence \(\gamma\) is a polynomial loop.

  To complete the proof of the proposition, we need to show that no
  fibre of \(P_\pol G \to G\) is empty, whence also no fibre of
  \(P_\per G \to G\) is empty.  As \(G\) is a compact, connected Lie
  group, for each \(g \in G\) there is some \(\xi \in \mf{g}\) such
  that \(\exp(\xi) = g\).  The path \(\eta_\xi\) lies in \(P_\pol G\)
  (and thus in \(P_\per G\)) and is in the fibre above \(g\).  Thus
  the fibres are non-empty.
\end{proof}

Proving that \(P_\per G\) is locally trivial is relatively
straightforward.  The case of \(P_\pol G\) is harder.  Therefore we
deal with \(P_\per G\) quickly now before passing to the -- for this
paper -- more relevant case of the polynomial loops in the next
section.

\begin{proposition}
  The space \(P_\per G\) is locally trivial over \(G\).
\end{proposition}

\begin{proof}
  To prove this, we require local sections.  Let \(g \in G\).  Let
  \(\xi \in \mf{g}\) be such that \(\exp(\xi) = g\).  Let \(\rho :
  [0,\frac12] \to [0, \frac12]\) be a smooth surjection which
  preserves the endpoints and is constant in a neighbourhood of each
  endpoint.  Let \(\phi : V \to U\) be a chart for \(G\) with \(U\) a
  neighbourhood of \(g\) such that \(\phi^{-1}(g) = 0\).

  For \(h \in U\), define a path \(\alpha_h : [0, 1] \to G\) by:
  \[
  \alpha_h(t) = \begin{cases}
    \exp(2\rho(t) \xi) & t \in [0, \frac12] \\
    \phi((2\rho(t - \frac12) + 1)\phi^{-1}(h)) & t \in [\frac12, 1]
	      \end{cases}
  \]
  By construction, \(\alpha_h\) is continuous.  Since \(\alpha_h\) is
  constant in a neighbourhood of \(\frac12\) and is smooth either
  side, it is smooth.  Moreover, as it is constant in neighbourhoods
  of \(0\) and \(1\), the concatenation \(\alpha_h \sharp (\alpha_h(1)
  \alpha_h)\) is smooth.  Hence \(\alpha_h\) extends via the formula:
  \[
  \alpha_h(t + n) = \alpha_h(1)^n \alpha_h(t)
  \]
  for \(t \in [0,1)\) and \(n \in \Z\), to a smooth path \(\R \to G\)
  such that \(\alpha_h(t + 1) = \alpha_h(1) \alpha_h(t)\) for all \(t
  \in \R\).

  Clearly, \(\alpha_h(1) = h\).  Also, the assignment \(h \to
  \alpha_h\) is smooth.  Therefore, \(h \to \alpha_h\) is a local
  section of \(P_\per G\) in a neighbourhood of \(g\).
\end{proof}

\subsubsection{The Polynomial Path Space}

The case of the polynomial path space is harder.  Regarding local
sections, it would appear from the definition that there are natural
local sections, namely \(g \to \eta_\xi\) where \(\exp(\xi) = g\).
However, except in the case of the unitary group, there is in general
no way to choose \(\xi\) smoothly in \(g\) for all points \(g \in G\)
(it is always possible to do so for an open dense subset, but this is
not good enough).

In fact, we are not able to prove that \(P_\pol G \to G\) is locally
trivial for all compact, connected \(G\) at this time.  The methods we
employ work on a case-by-case basis.  This is sufficient for our needs
as we are mainly interested in ordinary vector bundles with inner
products and thus in the structure groups \(U_n\) and \(S O_n\).  We
shall prove that \(P_\pol G \to G\) is locally trivial for these
groups and also for \(S U_n\).  There is no \emph{a priori} reason why
the argument for \(S O_n\) should not extend to \(S p_n\), using
quaternionic structures in place of complex structures but we feel
that this case is outside the focus of this paper.

The following result will prove useful in examining the structure of
\(P_\pol G\) in terms of \(P_\pol U_n\).

\begin{lemma}
  \label{lem:polsub}
  Let \(G\) be a compact, connected Lie group.  Consider \(G\) as a
  subgroup of \(U_n\).  Then \(P_\pol G = P_\per G \cap P_\pol U_n\).
\end{lemma}

\begin{proof}
  Clearly \(P_\pol G \subseteq P_\per G \cap P_\pol U_n\).  For the
  converse, let \(\alpha \in P_\per G \cap P_\pol U_n\).  Then
  \(\alpha = \eta_\xi \gamma\) for some \(\xi \in \mf{u}_n\) and
  \(\gamma \in L_\pol U_n\).  Now \(\exp(\xi) = \eta_\xi(1) =
  \alpha(1) \in G\) since \(\alpha \in P_\per G\).  Choose \(\zeta \in
  \mf{g}\) such that \(\exp(\zeta) = \exp(\xi)\).  Then:
  \[
  \alpha = \eta_\zeta \eta_{-\zeta} \eta_\xi \gamma.
  \]
  By proposition~\ref{prop:liepol}, \(\eta_{-\zeta} \eta_\xi\) is a
  polynomial loop in \(U_n\).  Since \(\alpha\) and \(\eta_\zeta\)
  both take values in \(G\), \(\eta_{-\zeta} \eta_\xi \gamma\) must
  also take values in \(G\).  It thus lies in \(L G \cap L_\pol U_n\)
  which is, by definition, \(L_\pol G\).  Therefore \(\alpha\) is of
  the form \(\eta_\zeta \beta\) with \(\zeta \in \mf{g}\) and \(\beta
  \in L_\pol G\).  Hence \(\alpha \in P_\pol G\).
\end{proof}

\subsubsection{The Unitary Group}

In the case of \(U_n\), there are local sections of the form \(g \to
\eta_\xi\) where \(\exp(\xi) = g\).  This will follow from
lemma~\ref{lem:comlie}.

\begin{proposition}
  \label{prop:unloctriv}
  The space \(P_\pol U_n\) is locally trivial over \(U_n\).
\end{proposition}

\begin{proof}
  Let \(s \in i \R\).  Let \(V_s \subseteq U_n\) be the open subset
  consisting of those operators which do not have \(-e^s\) as an
  eigenvalue.  Let \(\mf{v}_s \subseteq \mf{u}_n\) be the open subset
  consisting of those operators which have eigenvalue in the interval
  \((s - i \pi, s + i \pi)\).  The exponential map restricts to a
  diffeomorphism \(\exp : \mf{v}_s \to V_s\).  Let \(\log_s : V_s \to \mf{v}_s\)
  be its inverse.

  For a direct construction, define the \(s\)-logarithm \(\log_s : \T
  \ssetminus \{-e^s\} \to (s - i \pi, s + i \pi)\) as the inverse of
  the exponential map on this domain (note that this coincides with
  the above definition putting \(n = 1\)).  Let \(g \in V_s\).  Let
  \(E_1 \oplus \dotsb \oplus E_l\) be the orthogonal decomposition of
  \(\C^n\) into the eigenspaces of \(g\) with eigenvalues \(\lambda_1,
  \dotsc, \lambda_l\).  Then \(\log_s g\) is the operator which
  acts on \(E_j\) by multiplication by \(\log_s \lambda_j\).

  It is a simple exercise to show that \(\log_s g \in Z(g)\) for any
  \(g\) and \(s\) such that \(\log_s g\) is defined, that \(\log_s g\)
  is locally constant in \(s\), and that \(V_{s + 2 \pi i} = V_s\) and
  \(\mf{v}_{s + 2 \pi i} = \mf{v}_{s} + 2 \pi i 1_n\).

  The local sections of \(P_\pol U_n \to U_n\) are \(\alpha_s : V_s
  \to P_\pol U_n\) given by \(\alpha_s(g)(t) = \exp(t \log_s g)\).
\end{proof}

\subsubsection{The Special Unitary Group}

The method of the previous section works in \(U_n\) because every
point in \(U_n\) is \emph{exp-regular}; that is, is the image of a
point in \(\mf{u}_n\) such that the exponential map is a
diffeomorphism is a neighbourhood of that point.  This is not true for
a general Lie group.  It is straightforward to show that the preimage
of \(-1 \in S U_2\) under \(\exp : \mf{s u}_2 \to S U_2\) is a
countable number of copies of \(\CP^1\), hence \(-1 \in S U_2\) is not
exp-regular.

However, we can still prove that \(P_\pol S U_n \to S U_n\) is locally
trivial.  The strategy is to use the fact that there \emph{is} a point
in \(\mf{u}_n\) around which the exponential map is a local
diffeomorphism, and then use the fact that \(S U_n \to U_n \to S^1\)
is split.

\begin{proposition}
  The map \(P_\pol S U_n \to S U_n\) is locally trivial.
\end{proposition}

\begin{proof}
  Choose a unit vector \(v \in \C^n\).  Define the representation
  \(\sigma : \T \to U_n\) by \(\sigma(\lambda)v = \lambda v\) and
  \(\sigma(\lambda)\) is the identity on \(\langle v \rangle^\bot\).

  Let \(s \in i \R\).  Let \(V_s \subseteq U_n\) and \(\mf{v}_s \subseteq
  \mf{u}_n\) be as in the proof of proposition~\ref{prop:unloctriv}.
  Let \(\alpha_s : V_s \to P_\pol U_n\) be the local section defined
  in that proposition.

  Define \(\beta_s : V_s \cap S U_n \to P_\per U_n\) by:
  \[
  \beta_s(g)(t) = \alpha_s(g)(t) \,
  \sigma\Big(\det\big(\alpha_s(g)(-t)\big)\Big).
  \]

  Recall that \(\det \exp(\xi) = e^{\tr \xi}\).  Thus for \(g \in V_S
  \cap S U_n\):
  \[
  \det \alpha_s(g)(-t) = e^{\tr( -t \log_s(g))} = e^{-t \tr
  \log_s(g)}.
  \]
  As \(g \in S U_n\), \(e^{\tr \log_s(g)} = \det g = 1\) so \(\tr
  \log_s(g) = 2 \pi i k\) for some \(k \in \Z\).  Thus \(t \to e^{-t
  \tr \log_s(g)}\) is the map \(t \to z^{-k}\).  Hence \(\sigma( \det
  \alpha_s(g)(-t))\) is a polynomial loop in \(U_n\).  Thus
  \(\beta_s(g)(t) \in P_\pol U_n\).

  Then as \(\det \circ \sigma : \T \to \T\) is the identity, \(\det
  \beta_s(g)(t) = 1\) for all \(g, t\).  Hence \(\beta_s(g)(t) \in S
  U_n\) for all \(g, t\).  Thus by lemma~\ref{lem:polsub},
  \(\beta_s(g) \in P_\per S U_n \cap P_\pol U_n = P_\pol S U_n\).
\end{proof}

\subsubsection{The Special Orthogonal Group}
\label{sec:son}

The situation for \(S O_n\) is more complicated still.  The problem
here is with eigenvalue \(-1\).  It can be shown that \(g \in S O_n\)
is \emph{exp-regular} if and only if its \(-1\)-eigenspace has
dimension at most \(2\).  The solution comes from the theory of
\emph{unitary structures} which we now describe.

\begin{defn}
  Let \(E\) be a real vector space with an inner product.  A
  \emph{unitary structure} on \(E\) is an orthogonal map \(J : E \to
  E\) such that \(J^2 = -1\).
\end{defn}

\begin{proposition}
  \label{prop:cplxstr}
  Let \(E\) be a real even dimensional vector space with an inner
  product.  The properties of unitary structures that we shall need
  are:

  \begin{enumerate}
  \item \(E\) admits a unitary structure.

    \label{it:cplxeven}

  \item The set of unitary structures on \(E\) is \(O(E) \cap
    \mf{o}(E)\).

    \label{it:cplxorth}

  \item Let \(J\) be a unitary structure on \(E\).  Then \(\exp(\pi J)
    = -1_E\).

    \label{it:cplxexp}

  \item Let \(J_1, J_2\) be unitary structures on \(E\).
    Then:
    \(
    \eta_{-\pi J_1} \eta_{\pi J_2}
    \)
    is a polynomial loop in \(S O(E)\).

    \label{it:cplxpol}

  \item Let \(\xi \in \mf{s o}(E)\) be such that \(\xi\) does not have
    \(0\) as an eigenvalue.  Then there is a natural unitary structure
    \(J_\xi\) on \(E\) which varies smoothly in \(\xi\).  Considered
    as an element of \(\mf{s o}(E)\), \(J_\xi\) satisfies \([\xi,
    J_\xi] = 0\).  The assignment \(\xi \to J_\xi\) satisfies \(J_J =
    J\) (here \(J\) is considered as an element of \(\mf{s o}(E)\)),
    and \(J_{\xi + c J_\xi} = J_\xi\) for \(c > 0\).

    \label{it:cplxalg}

  \item Let \(g \in S O(E)\) be such that \(1\) is not an eigenvalue
    of \(g\).  Then \(\log_0(-g)\) is of the form \(\xi - \pi J_\xi\)
    for some \(\xi \in \mf{s o}(E)\) with \(\exp(\xi) = g\).

    \label{it:cplxdecomp}
  \end{enumerate}
\end{proposition}

In the last property we use the inclusion \(S O(E) \to U(E \otimes
\C)\) to define \(\log_0 : S O(E) \cap V_0 \to \mf{u}(E)\).  Since
\(\log_0\) commutes with complex conjugation%
\footnote{%
  It is the only one of the logarithms that we have defined with this
  property.%
}%
, the image of \(S O(E) \cap V_0\) lies in \(\mf{s o}(E)\).

\begin{proof}
  Property~\ref{it:cplxeven} is a standard property of complex
  structures whilst~\ref{it:cplxorth} is a simple deduction from the
  definition of a unitary structure.  Therefore we start with
  property~\ref{it:cplxexp}.

  \begin{enumerate}
    \addtocounter{enumi}{2}
  \item As an element of \(\mf{o}(E) = \mf{s o}(E)\), \(J\) is
    diagonalisable over \C.  Since \(J^2 = -1\), its eigenvalues are
    \(\pm i\).  Thus \(\pi J\) has eigenvalues \(\pm \pi i\).  Hence
    \(\exp(\pi J)\) has sole eigenvalue \(-1\).  As \(\exp(\pi J) \in
    S O(E)\), it is diagonalisable over \C and thus is \(-1_E\).

  \item This is a corollary of proposition~\ref{prop:liepol} together
    with the previous property.

  \item Diagonalise \(\xi\) over \C.  As \(\xi\) is a real operator,
    its eigenvalues and corresponding eigenspaces come in conjugate
    pairs.  As \(\xi\) is skew-adjoint, its eigenvalues lie on the
    imaginary axis in \C.  Let \(W \subseteq E \otimes \C\) be the sum
    of the eigenspaces of \(\xi\) corresponding to eigenvalues of the
    form \(i s\) with \(s > 0\).  Then \(\conj{W}\), resp.\
    \(W^\bot\), is the sum of the eigenspaces of \(\xi\) corresponding
    to eigenvalues of the form \(i s\) with \(s < 0\), resp.\ \(s \le
    0\).  The assumption on \(\xi\) implies that \(\conj{W} =
    W^\bot\).  Define \(J_\xi\) on \(E \otimes \C\) to be the operator
    with eigenspaces \(W\) and \(\conj{W}\) with respective
    eigenvalues \(i\) and \(-i\).  By construction, \(J^2 = -1\) and
    \(J^* J = 1\).  As the eigenspaces and eigenvalues of \(J\) come
    in conjugate pairs, \(J\) is a real operator and thus is a unitary
    structure.

    Since \(J_\xi\) is defined from the eigenspaces of \(\xi\), it
    varies smoothly in \(\xi\).  Moreover, as the eigenspaces of
    \(J_\xi\) decompose as eigenspaces of \(\xi\), \(J_\xi\) and
    \(\xi\) are simultaneously diagonalisable over \C.  Hence \([\xi,
    J_\xi] = 0\).

    It is clear from the construction that if \(\zeta\) and \(\xi\)
    can be simultaneously diagonalised and the eigenvalues of
    \(\zeta\) have the same parity on the imaginary axis as the
    corresponding ones of \(\xi\) then \(J_\zeta = J_\xi\).  In
    particular, \(J_J = J\) and \(J_{\xi + c J_\xi} = J_\xi\) for \(c
    > 0\).

  \item Let \(F\) be the \(-1\)-eigenspace of \(g\).  Then \(E\)
    decomposes \(g\)-invariantly as \(F \oplus F^\bot\).  As \(g\)
    does not have \(1\) as an eigenvalue, \(\log_0(-g)\) is
    well-defined.  Since the decomposition of \(E\) is
    \(-g\)-invariant:
    \[
    \log_0(-g) = \log_0(-g \restrict_{F}) + \log_0(-g
    \restrict_{F^\bot}) = \log_0(-g \restrict_{F^\bot}).
    \]
    This last step is because \(-g \restrict_F = 1_F\) so \(\log_0(-g
    \restrict_F) = 0_F\).

    Let \(\xi_{F^\bot} = \log_0(-g \restrict_{F^\bot})\).  As \(-g\)
    does not have \(1\) as an eigenvalue on \(F^\bot\),
    \(\xi_{F^\bot}\) does not have \(0\) as an eigenvalue.  Let
    \(J_{F^\bot}\) be the corresponding unitary structure.  As
    \([\xi_{F^\bot}, J_{F^\bot}] = 0\),
    \[
    \exp(\xi_{F^\bot} + \pi J_{F^\bot}) = \exp(\xi_{F^\bot}) \exp(
    \pi J_{F^\bot}) = (-g) \restrict_{F^\bot} (-1_{F^\bot}) = g
    \restrict_{F^\bot}.
    \]

    As \(g \in S O(E)\), \(F\) must be of even dimension.  Choose a
    unitary structure \(J_F\) on \(F\).  Then \(\exp(\pi J_F) = -1_F =
    g \restrict_F\).  Let \(\xi = \pi J_F + \xi_{F^\bot} + \pi
    J_{F^\bot}\).  Then:
    \[
    \exp(\xi) = \exp(\pi J_F) + \exp(\xi_{F^\bot} + \pi J_{F^\bot}) =
    -1_F + g \restrict_{F^\bot} = g.
    \]
    Then \(J_\xi = J_F + J_{F^\bot}\) so \(\xi - \pi J_\xi =
    \xi_{F^\bot}\), whence \(\xi - \pi J_\xi = \log_0 (-g)\).\qedhere
  \end{enumerate}
\end{proof}

\begin{theorem}
  \label{th:lctrso}
  The map \(P_\pol S O_n \to S O_n\) is locally trivial.
\end{theorem}

\begin{proof}
  We first describe a family of open sets which cover \(S O_n\).
  These will be the domains of the sections of \(P_\pol S O_n\).  The
  family is indexed by the interval \([-1,1]\) and by elements of \(S
  O_n\).

  Let \(r \in [-1,1]\).  Let \(W_r\) be the open subset of \(S O_n\)
  consisting of those \(g\) such that no eigenvalue of \(g\) (over \C)
  has real part \(r\).  For \(g \in W_r\) there is a \(g\)-invariant
  orthogonal decomposition of \(\R^n\) as \(E_{-1}^r(g) \oplus
  E_r^1(g)\) where the eigenvalues (over \C) of \(g\) on
  \(E_{-1}^r(g)\) have real part in the interval \([-1,r]\) and on
  \(E_r^1(g)\) in the interval \([r, 1]\).  Note that \(g\) cannot
  have eigenvalue \(1\) on \(E_{-1}^r(g)\), even if \(r = 1\), so as
  \(g \in S O_n\), \(E_{-1}^r(g)\) must have even dimension.

  Over each \(W_r\) is a vector bundle with fibre \(E_{-1}^r(g)\) at
  \(g\) (this will have different dimension on the different
  components of \(W_r\)).  Over most \(W_r\)'s this bundle is not
  trivial.  Therefore we find smaller open sets over which we can
  trivialise it.

  Let \(r \in [-1,1]\) and \(g \in W_r\).  Define \(W_r(g)\) to be the
  open subset of \(S O_n\) consisting of those \(h \in W_r\) for which
  the orthogonal projection \(E_{-1}^r(h) \to E_{-1}^r(g)\) is an
  isomorphism.

  Over \(W_r(g)\), therefore, the aforementioned vector bundle is
  trivial and of constant even dimension.  Hence, we can choose a
  unitary structure \(J_h\) on each \(E_{-1}^r(h)\) which varies
  smoothly in \(h\).

  Extend \(J_h\) to a skew-adjoint operator on \(\R^n\) by defining it
  to be zero on \(E_r^1(h)\).  Let \(\epsilon(h) = h \exp(-\pi J_h)
  \in S O_n\).  Then \(\epsilon(h)\) agrees with \(h\) on \(E_r^1(h)\)
  and is \(-h\) on \(E_{-1}^r(h)\).  Since \(h\) does not have
  eigenvalue \(-1\) on \(E_r^1(h)\) and does not have eigenvalue \(1\)
  on \(E_{-1}^r(h)\), \(\epsilon(h)\) does not have eigenvalue \(-1\)
  on \(\R^n\) and so lies in the domain of \(\log_0\).  Also, as
  \(J_h\) varies smoothly in \(h\), \(h \to \epsilon(h)\) is smooth.

  Define \(\beta_{r,g} : W_r(g) \to P_\per S O_n\) by:
  \[
  \beta_{r,g}(h)(t) = \exp\big( t \log_0( \epsilon(h)) \big) \exp(t
  \pi J_h).
  \]
  This is a smooth path in \(S O_n\) since both \(\log_0 (
  \epsilon(h))\) and \(J_h\) lie in \(\mf{s o}_n\).  It varies
  smoothly in \(h\) since both \(\epsilon(h)\) and \(J_h\) are smooth
  in \(h\).  Since \(\epsilon(h) = h \exp(-\pi J_h)\),
  \(\beta_{r,g}(h)(1) = h\) so it is a path above \(h\).  We need to
  show that it lies in \(P_\pol S O_n\).

  Now \(\epsilon(h)\) respects the decomposition \(E_{-1}^r(h) \oplus
  E_r^1(h)\) of \(\R^n\), therefore so does \(\log_0(\epsilon(h))\).
  Accordingly, write \(\log_0(\epsilon(h)) = \xi_{-1}^r + \xi_r^1\).

  Consider the situation on \(E_{-1}^r(h)\).  Since \(\exp(\xi_{-1}^r)
  = \epsilon(h) = -h\) (all restricted to \(E_{-1}^r(h)\)), by
  property~\ref{it:cplxdecomp}, \(\xi_{-1}^r = \zeta - \pi J_\zeta\)
  for some \(\zeta \in \mf{s o}(E_{-1}^r(h))\) with \(\exp(\zeta) =
  h\).
  Extend \(J_\zeta\) to \(\R^n\) by defining it to be zero on
  \(E_r^1(h)\).  Let \(\xi = \zeta + \xi_r^1\).  Then \(\exp(\xi) =
  h\), \([\xi, J_\zeta] = 0\), and \(\log_0( \epsilon(h)) = \xi - \pi
  J_\zeta\).  Therefore:
  \[
  \beta_{r,g}(h)(t) = \exp(t \xi) \exp(-t \pi J_\zeta) \exp(t \pi
  J_h).
  \]

  Since \(J_\zeta\) and \(J_h\) are both extensions to \(\R^n\) by
  zero of unitary structures on \(E_{-1}^r(h)\), then by
  property~\ref{it:cplxpol},  \(\exp(-t \pi J_\zeta) \exp(t \pi
  J_h)\) is a polynomial loop in \(S O_n\).  Hence \(\beta_{r,g}(h)\)
  lies in \(P_\pol S O_n\).
\end{proof}

\subsection{The Polynomial Vector Bundles}
\label{sec:polvect}

Now that we have principal bundles, given a representation we can
construct vector bundles.  Let \(V\) be a finite dimensional vector
space with an inner product, either real or complex.  Let \(L V\) be
the space of smooth loops in \(V\) and \(L_\pol V\) the space of
polynomial loops.  If \(V\) is complex then \(L_\pol V = V
[z^{-1},z]\); if \(V\) is real then \(L_\pol V = L V \cap L_\pol (V
\otimes \C)\).

Let \(G\) be a compact, connected Lie group which acts on \(V\) by
isometries.  In the polynomial case, assume that \(G\) is one of
\(U_n\), \(S U_n\), or \(S O_n\).  Then \(L G\) acts on \(L V\) and
\(L_\pol G\) acts on \(L_\pol V\).  Therefore we have vector bundles
over \(G\) together with a bundle inclusion:
\[
P_\pol V := P_\pol G \times_{L_\pol G} L_\pol V \to
P_\per V := P_\per G \times_{L G} L V.
\]

We shall now give an alternative view of these vector bundles which
will be more enlightening in terms of their structure.

Let \(P V\) be the full path space of \(V\).  Define \(\tau : P V \to
P V\) to be the shift operator: \((\tau \beta )(t) = \beta(t + 1)\).
Let \(D\) denote the differential operator: \((D \beta)(t) =
\diff{\beta}{t}(t)\).  There is a strong connection between these
operators: \(D\) is the infinitesimal generator of the group of
translations on \(P V\) and \(\exp(D) = \tau\).

The motivation for considering these operators is that they give
simple descriptions of \(L V\) and \(L_\pol V\) inside \(P V\).  The
loop space, \(L V\), is the \(+1\)-eigenspace of \(\tau\).   The space
of polynomial loops inside \(L V\) is the union of the finite
dimensional \(D\)-invariant subspaces of \(L V\).

In the complex case, we can write this as the linear span of the
eigenvectors of \(D\).  This does not carry over to the real case,
however, as the only eigenvectors of \(D\) are the constant maps.

\begin{theorem}
  Let \(g\) in \(G\).  The fibre of \(P_\per V\) above \(g\) is the
  space of \(\phi \in P V\) such that \(\tau \phi = g \phi\).

  The fibre of \(P_\pol V\) above \(g\) is the union of the finite
  dimensional \(D\)-invariant subspaces of the fibre of \(P_\per V\)
  above \(g\).
\end{theorem}

\begin{proof}
  An element of \(P_\per V\) in the fibre above \(g\) is represented
  by a pair \((\alpha, \beta)\) with \(\alpha \in P_\per G\) above
  \(g\) and \(\beta \in L V\).  Any alternative representative is of
  the form \((\alpha \gamma, \gamma^{-1} \beta)\) for some \(\gamma
  \in L G\).

  Thus the map \(\phi : \R \to V\) defined by \(\phi := \alpha \beta\)
  depends only on the element of \(P_\per V\) and not on the choice of
  representative.  This satisfies:
  \[
  (\tau \phi)(t) = \phi(t + 1) = \alpha(t + 1) \beta(t + 1) = g
  \alpha(t) \beta(t) = g \phi(t).
  \]
  Hence \(\tau \phi = g\phi\).

  Conversely, suppose that \(\tau \phi = g \phi\).  Choose some
  \(\alpha \in P_\per G\) above \(g\) and define \(\beta :=
  \alpha^{-1} \phi\).  Then \(\beta(t + 1) = \alpha^{-1}(t) g^{-1} g
  \phi(t) = \beta(t)\) so \(\beta \in L V\).  Changing \(\alpha\) to
  \(\alpha \gamma\) changes \(\beta\) to \(\gamma^{-1} \beta\).  Hence
  the element in \(P_\per V\) represented by \((\alpha, \beta)\)
  depends only on \(\phi\).

  Now we consider the polynomial path space.  We need to show that the
  fibre of \(P_\pol V\) above \(g\) is the union of the finite
  dimensional subspaces of the fibre of \(P_\per V\) that are
  \(D\)-invariant.

  Let \(\xi \in \mf{g}\) be such that \(\exp(\xi) = g\).  We consider
  two actions of \(\xi\) on \(P V\).  The first is the isomorphism
  \(\alpha \to \eta_{-\xi} \alpha\) which maps \(P_{\per,g} V\) onto
  \(L V\).  The second is \(\alpha \to \xi \alpha\), extending the
  action of \(\mf{g}\) on \(V\) to \(P V\).  As \(\xi\) is a finite
  dimensional operator, it has a minimum polynomial.  This is true
  also of its action on \(P V\).  Therefore any finite dimensional
  subspace of \(P V\) is contained in a finite dimensional
  \(\xi\)-invariant subspace.  Moreover, the action of \(\xi\) on \(P
  V\) commutes with that of \(D\) so any finite dimensional
  \(D\)-invariant subspace of \(P V\) is contained in a finite
  dimensional subspace that is both \(D\)-invariant and
  \(\xi\)-invariant.

  Hence as \(\xi\) preserves both \(P_{\per, g} V\) and \(L V\), when
  considering the union of finite dimensional \(D\)-invariant
  subspaces in either, it is sufficient to consider those that are in
  addition \(\xi\)-invariant.

  We shall now show that \(W \subseteq L V\) is \(\xi\) and
  \(D\)-invariant if and only if \(\eta_\xi W\) is \(\xi\) and
  \(D\)-invariant.  This will establish the result.

  The \(\xi\)-invariance is straightforward since \(\xi\) commutes
  with \(\eta_\xi\).  Hence \(W \subseteq L V\) is \(\xi\)-invariant
  if and only if \(\eta_\xi W \subseteq P_{\per,g} V\) is
  \(\xi\)-invariant.

  If \(W\) is \(\xi\) and \(D\)-invariant, then consider \(\alpha \in
  \eta_{\pm \xi} W\) (the \(\pm\) allows us to consider both
  directions at once).  This is of the form \(\eta_{\pm \xi} \beta\)
  for some \(\beta \in W\).  Then:
  \[
  D \alpha = (D \eta_{\pm \xi}) \beta + \eta_{\pm \xi} (D \beta) =
  \eta_{\pm \xi} ( \pm \xi \beta + D \beta) \in \eta_{\pm \xi} W.
  \]
  Hence \(\eta_{\pm \xi} W\) is \(D\)-invariant.
\end{proof}

An immediate corollary of this is that the fibres of \(P_\per V\) and
of \(P_\pol V\) are \(D\)-invariant.  For \(P_\per V\) this follows
from the fact that \(\exp(D) = \tau\) so \(D\) and \(\tau\) commute.
If we wish to emphasise the fibre, we shall refer to \(D\) as \(D_g\).

In the complex case, as \(D_g\) is skew-adjoint, any element of the
fibre of \(P_\pol V\) above \(g\) is thus the sum of eigenvectors of
\(D_g\).

When viewing a fibre of \(P_\per V\) or \(P_\pol V\) as a subspace of
\(P V\), the corresponding element \(g \in G\) is not uniquely
determined by any one path (contrast with the case of \(P_\per G\) or
\(P_\pol G\)).  Thus to keep track of the fibre, we shall often use
the notation \((g, \phi)\).

There is an action of \(G\) on \(P_\per V\) and on \(P_\pol V\) given
by the following equivalent definitions:
\begin{align*}
g \cdot [\alpha, \beta] &= [g \alpha, \beta], \\
g \cdot [\alpha, \beta] &= [g \alpha g^{-1}, g \beta], \\
g \cdot (h, \phi) &= (g h g^{-1}, g \phi).
\end{align*}
We put in both of the top two descriptions to show that the two
actions of \(G\) on \(P_\per G\) (and thus on \(P_\pol G\)) define the
same action on \(P_\per V\) (and \(P_\pol V\)).  This action preserves
the sub-bundle \(P_\pol V\) and sends the operator \(D_h\) to \(D_{g h
g^{-1}}\).

\newpage
\section{Dual Loop Bundles}
\label{sec:dual}

The goal of this section is to construct the inner product and the
Hilbert completion of the dual of the vector bundle \(L E \to L M\),
where \(E \to M\) is a real or complex vector bundle.  The first part
of this construction involves defining the polynomial loop bundle,
\(L_\pol E \to L M\) and proving that it is a locally trivial vector
bundle modelled on \(L_\pol \F^n\), for \F one of \R or \C.  Once this
has been defined, we thicken it to a Hilbert bundle which is a
sub-bundle of \(L E\).  This dualises to the required completion of
\(L^* E\).  We show how to construct an inner product on this bundle
by finding an isomorphism of the completion of \(L E\) with the
completion of \(L^* E\).

In section~\ref{sec:polprop} we discuss the basic properties of the
polynomial loop bundle, and thus of the Hilbert completion of \(L^*
E\).  In particular we consider the action of the group of
diffeomorphisms of the circle.  The natural action on \(L E\) does not
preserve the polynomial sub-bundle but it can be modified to an action
which does.

The construction of the polynomial loop bundle relies on the holonomy
operator coming from a connection.  The holonomy map can be viewed as
a variant of a classifying map for the original loop bundle.  In
section~\ref{sec:twist} we examine this idea.

\subsection{Polynomial Loop Bundles}
\label{sec:poly}

Let \(M\) be a smooth finite dimensional manifold without boundary.
Let \(G\) be one of \(U_n\), \(S U_n\), or \(S O_n\).  Let \F be the
corresponding field.  Let \(Q \to M\) be a principal \(G\)-bundle.
Let \(E = Q \times_G \F^n\) be the corresponding vector bundle.  As
\(G\) preserves the inner product on \(\F^n\), \(E\) carries a
fibrewise inner product.  Let \(\nabla\) be a covariant differential
operator on \(E\) coming from a connection on \(Q\).

We think of a point in a fibre \(Q_p\) as being an isometry \(\F^n \to
E_p\).  We shall also use the \emph{adjoint bundle} associated to
\(Q\), \(Q^\ad := Q \times_{\text{conj}} G\) where \(G\) acts on
itself by conjugation.  This is a bundle of groups.  A point in a
fibre \(Q^\ad_p\) is an isometry of \(E_p\) to itself.

It is a standard result that the loop and path spaces of \(E\) form
vector bundles over, respectively, the loop and path spaces of \(M\)
with frame bundles the loop and path spaces of \(Q\) and adjoint
bundles the loop and path spaces of \(Q^\ad\).

Recall from section~\ref{sec:notepath} that for \(X\) each of \(E\),
\(Q\), and \(Q^\ad\), the fibre of the bundle \(P^M X\) above \(p \in
M\) is \(P (X_p)\).  Thus:
\begin{align*}
P^M E &= E \otimes \Ci(\R,\F), & L^M E &= E
\otimes L \F, \\
&= Q \times_G \Ci(\R, \F^n), && = Q \times_G L \F^n,\\
P^M Q &= Q \times_{G} P G, & L^M Q &= Q \times_{G} L G,\\
P^M Q^\ad &= Q \times_{\text{conj}} P G, & L^M Q^\ad &= Q
\times_{\text{conj}} L G.
\end{align*}
As with \(M\) inside \(L M\) and \(P M\), \(G\) sits inside \(L
G\) and \(P G\) as the constant loops.  In the middle line, the
action is as a subgroup, in the third line the action is via
conjugation.

Since \(G\) acts on \(P_\per G\), \(P_\pol G\), \(P_\per \F^n\), and
\(P_\pol \F^n\), we can define corresponding bundles over \(M\).  For
\ttype each of \tper or \tpol, let:
\begin{align*}
  P^M_\type E &:= Q \times_G P_\type \F^n, \\
  P^M_\type Q &:= Q \times_G P_\type G, \\
  P^M_\type Q^\ad &:= Q \times_{\text{conj}} P_\type G.
\end{align*}
In the middle line, we use the left action of \(G\) on \(P_\type G\).
In the last line, we use the conjugation action.  As each of the model
spaces for these bundles is itself a bundle over \(G\) and the actions
on the total spaces induce the conjugation action on the base, for
\(X\) each of \(E\), \(Q\), and \(Q^\ad\), \(P^M_\type X\) is a fibre
bundle over \(Q^\ad\) with fibre \(L_\type Y\), where \(Y\) is either
\(\F^n\) or \(G\) as appropriate.

The covariant differential operator defines a parallel transport
operator.  This defines three compatible families of bundle maps
\(\psi_X^t : X^t \to P X\), for \(X\) each of \(E\), \(Q\), and
\(Q^\ad\).  The properties of these maps are:
\begin{align}
\psi_E^t(p q w) &= \psi_{Q^\ad}^t(p) \psi_Q^t(q) w, && p \in
Q^{\ad,t}_\gamma,\, q \in Q^t_\gamma,\, w \in \F^n \subseteq P
\F^n. 
\label{eq:psitriple}
\\
\psi_X^t e_t \psi_X^s &= \psi_X^s,
\label{eq:psibreak}
\\
e_{s + 1} \psi_X^{t + 1} &= e_s \psi_X^t, && \text{over } L M.
\label{eq:psiloop}
\end{align}
For the second, note that \(e_t \psi_X^s\) is a map from \(X^s\) to
\(X^t\).  This compatibility relation is the statement that if one
parallel transports from time \(s\) to time \(t\) and then on from
time \(t\) to some-when else, it is the same as transporting straight
from \(s\) to ones final time.  For the last, over \(L M\) then \(X^{t
+ 1} = X^t\) so the domains and codomains of these maps are the same.
This property is then an application of the fact that the parallel
transport operator is intrinsic to \(M\), therefore the parallel
transport from \(X^t\) to \(X^s\) is the same as that from \(X^{t +
1}\) to \(X^{s + 1}\).

These operators extend to bundle equivalences:
\begin{equation}
\label{eq:buneq}
\Psi_X^t : P^{M,t} X \to P X, 
\end{equation}
with the property that \(e_s( \Psi_X^t \alpha) = (e_s \psi_X^t)
(\alpha(s))\).  Note that these equivalences have been chosen such
that \((\Psi_X^t \alpha)(s)\) always lies in \(X^s\) no matter which
\(t\) was the starting point.

Recall that, for \(X\) each of \(E\), \(Q\), or \(Q^\ad\), \(L X\)
sits inside \(P^L X\).  It is straightforward to recognise this
submanifold: \(L X\) consists of those paths \(\beta \in P^L X\) which
are themselves periodic.  Note that for any path \(\beta\) in \(P^L
X\) then \(\beta(t + 1)\) and \(\beta(t)\) both lie in the same fibre
of \(X \to M\).

Thus in the right-hand side of~\eqref{eq:buneq} (restricted to \(L
M\)), it is straightforward to recognise the sub-bundles consisting of
the loops.  We wish to transfer this recognition principle to the
left-hand side of~\eqref{eq:buneq}.  We do this using the
\emph{holonomy operator}.

\begin{defn}
  On \(L M\), define the fibrewise operators \(h_X : X^0 \to X^0\) by
  \(h_X = e_1 \psi_X^0\).
\end{defn}

Over \(P M\), \(e_1 \psi_X^0\) is a map \(X^0 \to X^1\).  Over \(L M\)
then \(X^0 = X^1\) so \(h_X\) is as defined.  The fibres of \(Q^\ad\)
act on each of \(E\), \(Q\), and \(Q^\ad\): on \(E\) the action is by
definition, on \(Q\) and on \(Q^\ad\) by composition.

\begin{lemma}
  The operator \(h_E\) is a section of \(Q^{\ad,0}\).  The operators
  \(h_E\), \(h_Q\), and \(h_{Q^\ad}\) satisfy: \(h_{Q^\ad}(p) h_E =
  h_E p\), and \(h_Q(q) = h_E q\).  Thus \(h_E\) determines both
  \(h_Q\) and \(h_{Q^\ad}\).
\end{lemma}

\begin{proof}
  Since \(e_1 \psi^0_E\) is a fibrewise isometry
  \(E^0 \to E^0\), it is a section of \(Q^{\ad,0}\).  Then
  from~\eqref{eq:psitriple}, for \(p \in Q^{\ad,0}\), \(q \in Q^0\),
  \(v \in E^0\), and \(w \in \F^n \subseteq P \F^n\):
  \begin{align*}
    (h_E p)v &= (e_1 \psi^0_E  p ) v \\
    &= e_1 (\psi^0_E (p v)) \\
    &= e_1 (\psi^0_{Q^\ad}(p) \psi^0_E(v)) && \text{by
    \nmeqref{eq:psitriple}} \\
    &= (e_1 \psi^0_{Q^\ad})(p) (e_1 \psi^0_E)(v) \\
    &= h_{Q^\ad}(p) h_E(v). \\
    (h_E q)w &= (e_1 \psi^0_E q) w \\
    &= e_1( \psi^0_E( q w)) \\
    &= e_1( \psi^0_Q(q) w ) && \text{by \nmeqref{eq:psitriple}}\\
    &= (e_1 \psi^0_Q)(q) w \\
    &= h_Q(q) w. \qedhere
  \end{align*}
\end{proof}

\begin{lemma}
  \(e_{t + 1} \psi_X^0 = e_t \psi^0_X h_X\).
\end{lemma}

\begin{proof}%
  \begin{align*}
    e_{t + 1} \psi_X^0 &= e_{t + 1} \psi^1_X e_1 \psi^0_X && \text{by 
    \nmeqref{eq:psibreak}} \\
    &= e_t \psi^0_X e_1 \psi^0_X && \text{by \nmeqref{eq:psiloop}} \\
    &= e_t \psi^0_X h_X. \qedhere
  \end{align*}
\end{proof}

\begin{corollary}
  Under the bundle isomorphism of~\eqref{eq:buneq}, the sub-bundle \(L
  X\) of \(P^L X\) corresponds to:
  \[
  \{\alpha(t) \in P^{M,0} X : h_X \alpha(t + 1) = \alpha(t)\}.
  \]
\end{corollary}

\begin{proof}
  An element \(\alpha \in P^{M,0} X\) is mapped to a loop in \(P^L X\)
  if and only if \((\Psi^0_X \alpha)(t + 1) = (\Psi^0_X \alpha)(t)\)
  for all \(t \in \R\).  The left-hand side of this simplifies to:
  \[
  e_{t + 1} (\Psi^0_X \alpha) = (e_{t + 1} \psi^0_X)(\alpha(t + 1)) =
  (e_t \psi^0_X)(h_X \alpha(t + 1))
  \]
  whilst the right-hand side simplifies to:
  \[
  e_t (\Psi^0_X \alpha) = (e_t \psi^0_X)(\alpha(t)).
  \]

  Since \(e_t \psi^0_X : X^0 \to X^t\) is an isomorphism, this implies
  that \(\Psi^0_X \alpha\) is a loop if and only if \(h_X \alpha(t +
  1) = \alpha(t)\) for all \(t \in \R\).
\end{proof}

A section \(\chi : L M \to Q^{\ad,0}\) is the same thing as a map
\(\chi : L M \to Q^\ad\) which covers the map \(e_0 : L M \to M\).
For such a section, \(X\) each of \(Q\), \(Q^\ad\), or \(E\), and
\ttype each of \tpol or \tper, let \(L^{M,0,\chi}_\type X \to L M\) be
the pull-back of \(P^M_\type X \to Q^\ad\) via the map \(\chi : L M
\to Q^\ad\).

\begin{corollary}
  For \(X\) each of \(Q\), \(Q^\ad\), and \(E\), \(\Psi^0_X\)
  restricts to a bundle isomorphism \(L^{M,0,h_E^{-1}}_\per X \to L
  X\).
\end{corollary}

\begin{defn}
  The \emph{polynomial loop bundles}, \(L_\pol X\), for \(X\) each of
  \(Q\), \(Q^\ad\), and \(E\) are defined to be the images in \(L X\)
  of \(L^{M,0,h_E^{-1}}_\pol X\) under the map \(\Psi^0_X\).
\end{defn}

The following is immediate:

\begin{proposition}
  The polynomial loop bundles are locally trivial with \(L_\pol Q\) a
  \(L_\pol G\)-principal bundle, \(L_\pol Q^\ad\) a bundle of groups
  modelled on \(L_\pol G\), and \(L_\pol E\) a vector bundle modelled
  on \(L_\pol \F^n\).  Moreover:
  \begin{align*}
  L_\pol Q^\ad &= L_\pol Q \times_{\text{conj}} L_\pol G, \\
  L_\pol E &= L_\pol Q \times_{L_\pol G} L_\pol \F^n, \\
  L Q &= L_\pol Q \times_{L_\pol G} L G, \\
  L Q^\ad &= L_\pol Q \times_{\text{conj}} L G, \\
  L E &= L_\pol Q \times_{L_\pol G} L \F^n.
  \end{align*}
\end{proposition}

The bundle \(L_\pol E\) has a more concrete description in terms of
the connection on \(E\).  For any path \(\gamma : \R \to M\), the
connection on \(E\) defines a covariant differential operator
\(D_\gamma : \Gamma_\R(\gamma^* E) \to \Gamma_\R(\gamma^* E)\); that
is, \(D_\gamma : P_\gamma E \to P_\gamma E\).  As the map \(\Psi_E^0\)
was constructed using parallel transport, it (rather, its inverse)
takes \(D_\gamma\) to the operator \(\diff{}{t}\) acting on \(P^{M,0}
E\).  If \(\gamma\) happens to be a loop, \(D_\gamma\) restricts to an
operator on \(L_\gamma E\).  As \(\Psi_E^0\) identifies \(L_\gamma E\)
with the fibre of \(P_\per E \to Q^{\ad,0}\) above
\(h_E^{-1}(\gamma)\), it takes \(D_\gamma\) to the operator
\(D_{h_E^{-1}(\gamma)}\).

Hence \(L_{\pol, \gamma} E\) can be constructed from the action of
\(D_\gamma\) on \(L_\gamma E\) in the same fashion as \(P_{\pol,g}
\F^n\) from \(P_{\per,g} \F^n\), namely as the union of the finite
dimensional \(D_\gamma\)-invariant subspaces of \(L_\gamma E\).  In
the complex case, \(L_{\pol, \gamma} E\) is the span of the
eigenvalues in \(L_\gamma E\) of \(D_\gamma\).

\subsection{The Completion of the Cotangent Bundle}
\label{sec:cotangent}

The Hilbert completion of the cotangent bundle is now straightforward.
We merely need to select a Hilbert space that lies between \(L_\pol
\R^n\) and \(L \R^n\).  The group \(L_\pol S O_n\) will act on this
Hilbert space and thus we can construct a locally trivial bundle over
the loop space \(L M\) with fibre a Hilbert space which sits naturally
between \(L_\pol T M\) and \(L T M\).  Dualising this will yield a
Hilbert space sitting naturally between \(L^* T M\) and \(L_\pol^* T
M\).  On fibres, this will be a Hilbert completion of \(T^* L M =
L^* T M\).

There are many Hilbert spaces between \(L_\pol \R^n\) and \(L \R^n\).
We choose \(L^2_e \R^n\).  The choice of \(e\) is dictated by the
desire not to have unnecessary constants at a later stage.  It is not
overly significant.  There is an obvious inner product on \(L^2_e
\R^n\) but as it is not preserved under the action of \(L_\pol S
O_n\), we shall not spend any time discussing it.

\begin{lemma}
  Let \(G\) be one of \(S O_n\), \(S U_n\), or \(U_n\).  Let \F be \R
  or \C as appropriate.  Then \(L_\pol G\) acts continuously on
  \(L^2_e \F^n\).
\end{lemma}

\begin{proof}
  It is sufficient to show that \(L_\pol U_n\) acts on \(L^2_e \C^n\)
  since \(L_\pol S U_n\) and \(L_\pol S O_n\) are subgroups of
  \(L_\pol U_n\) and the action of \(L_\pol S O_n\) on \(L^2_e \R^n\)
  is well-defined if and only if its action on \(L^2_e \C^n\) is
  well-defined.

  To show that \(L_\pol U_n\) acts on \(L^2_e \C^n\), it is sufficient
  to show that \(L_\pol M_n(\C)\) acts.  This is straightforward since
  it is generated as an algebra by \(M_n(\C)\) and \(z\) which both
  act on \(L^2_e \C^n\).
\end{proof}

Thus given a vector bundle \(E \to M\) with structure bundle \(Q\) we
obtain a Hilbert bundle \(L^2_e E \to L M\) as \(L_\pol Q
\times_{L_\pol G} L^2_e \F^n\).  There are maps \(L_\pol E \to L^2_e E
\to L E\) which locally look like \(L_\pol \F^n \to L^2_e \F^n \to L
\F^n\).  Dualising this bundle yields the required completion of \(L^*
E\).

As remarked above, the natural inner product on \(L^2_e \F^n\) is not
preserved by the action of \(L_\pol G\).  It is a simple matter to
show that when the circle action is taken into account, the only
Hilbert completion of \(L_\pol \F^n\) on which \(L_\pol G\) can act
unitarily is \(L^2 \F^n\).  Therefore, we shall have to find another
route to an inner product on \(L^2_e E\) (and thus its dual).  The
route we choose is to construct an isomorphism of Hilbert bundles
\(L^2_e E \to L^2 E\).  This will allow us to pull-back the inner
product on \(L^2 E\) to \(L^2_e E\).

\begin{proposition}
  There is a well-defined bundle map \(L_\pol E \to L_\pol E\) given
  by:
  \[
  \alpha \to (\cos D_\gamma) \alpha
  \]
  where:
  \[
  \cos D_\gamma = \sum_{j = 0}^\infty \frac1{(2 j)!} D_\gamma^{2 j}.
  \]

  This extends to an isomorphism of Hilbert bundles \(L^2_e E \to L^2
  E\).
\end{proposition}

\begin{proof}
  We start with the fibrewise situation.  The fibre of \(L_\pol E\)
  above a loop \(\gamma\) is the union of \(D_\gamma\)-invariant
  finite dimensional subspaces of \(L_\gamma E\).  On any finite
  dimensional space, the power series denoted by \(\cos A\) converges
  for any operator \(A\).  Therefore, \(\cos D_\gamma\) is
  well-defined on each finite dimensional \(D_\gamma\)-invariant
  subspace of \(L_\gamma E\) and hence on \(L_\pol E\).

  When considering the Hilbert completions of \(L_{\pol, \gamma} E\),
  it is sufficient to assume that \(E\) is complex.  In this case,
  there is a basis for \(L_{\pol, \gamma} E\) of eigenvectors of
  \(D_\gamma\).  We can choose this basis to have the following
  properties:
  \begin{enumerate}
  \item There are \(n\) eigenvectors \(v_1, \dotsc, v_n \in L_\gamma
    E\) such that the corresponding eigenvalues are of the form \(i
    s_1, \dotsc, i s_n\) with each \(s_j \in [0, 2 \pi i)\).

  \item The other eigenvectors are of the form \(z^k v_j\) for some
    \(k \in \Z\).
  \end{enumerate}
  The eigenvalue of \(z^k v_j\) is \(i s_j + 2 \pi i k\).  Therefore,
  \((\cos D_\gamma) z^k v_j = \cosh(s_j + 2 \pi k) z^k v_j\).

  We wish to describe the sequences \((a^j_k)_{k \in \Z, j = 1,
  \dotsc, n}\) such that \((a^j_k \cosh(s_j + 2 \pi k))\) is
  square-summable.  It is sufficient to consider each \(j\) in turn so
  we consider \((a_k)_{k \in \Z}\) such that \((a_k \cosh(s + 2 \pi
  k))\) is square-summable for some \(s \in [0, 2\pi)\).

  Now elementary analysis shows that for all \(x \in \R\):
  \[
  \cosh (s) \ge \frac{\cosh (x + s)}{e^{\abs{x}}} \ge \frac12
  \min\{e^s, e^{-s}\}.
  \]
  Therefore, \((a_k \cosh(s + 2 \pi k))\) is square-summable if and
  only if \((a_k e^{2 \pi \abs{k}})\) is square-summable.  This is
  precisely the condition that the loop corresponding to \((a_k)\)
  extend analytically over an annulus of radii \(e\) and \(e^{-1}\)
  and be square-summable on the boundaries.
\end{proof}

Therefore, although we cannot transfer the standard inner product on
\(L^2_e \F^n\) to the fibres of \(L^2_e E\), we can use the operator
\(D_\gamma\) to define an inner product on the fibres of \(L^2_e E\)
which is equivalent to the standard inner product on \(L^2_e \F^n\).
This inner product is defined by:
\[
\ipv{\alpha}{\beta} = \ip{(\cos D_\gamma) \alpha}{(\cos D_\gamma)
  \beta}
\]
where \(\ipc\) is the inner product on \(L^2_\gamma E\).  Thus
\(\alpha \to (\cos D_\gamma) \alpha\) is an isometry from
\(L^2_{e,\gamma} E\) to \(L^2_\gamma E\).

Since an inner product on a Hilbert space defines a (conjugate-linear)
isomorphism with its dual, we can extend the inclusion \(L^* E \to
L^{2,*}_e E\) to a triple \(L^* E \to L^{2,*}_e E \cong L^2_e E \to L
E\).  This mirrors the standard triple \(L E \to L^2 E \cong L^{2,*} E
\to L^* E\) coming from the standard inner product on \(L E\).
Putting these together (and suppressing the isomorphisms) yields a
non-commuting square:
\[
\begin{CD}
  L E @>>> L^2 E \\
  @AAA @VVV \\
  L^2_e E @<<< L^* E .
\end{CD}
\]

As mentioned above, the group \(L_\pol G\) does not act on \(L^2_e
\F^n\) by isometries but does so act on \(L^2 \F^n\).  Therefore the
isomorphism \(L^2_e E \to L^2 E\) can be viewed as a careful
alteration of the structure group of \(L^2_e E\) to one which acts by
isometries.  We can do the same with \(L^* E\).

We have constructed a chain of bundle maps which on fibres looks like:
\[
L^2_e \F^n \to L \F^n \to L^2 \F^n \cong L^{2,*} \F^n \to L^* \F^n \to
L^{2,*}_e \F^n.
\]
The isometry \(L^2_e E \cong L^2 E\) takes this chain to one which
ends in \(L^{2,*} E\).  On fibres this chain is obtained from the one
above by adding a subscript \(e\) to each space.  Using the identities
\((L_e^2)_e \F^n = L_{e^2}^2 \F^n\) and \((L_e^{2,*})_e \F^n = L^{2,*}
\F^n\), the fibres of the vector bundles in this new chain are:
\[
L^2_{e^2} \F^n \to L_e \F^n \to L^2_e \F^n \cong L^{(2,*)}_e \F^n
\to L^*_e \F^n \to L^{2,*} \F^n.
\]

In this picture, we have identified \(L^* E\) with the sub-bundle
\(L^*_e E\) of \(L E\) and then taken the inner product and completion
of \(L^* E\) to be \(L^{2,*} E\).  Since the action of \(L_\pol G\) on
\(L^*_e E\) is by isometries, this identification of \(L^* E\) with
\(L^*_e E\) can be viewed as a careful alteration of the structure
group of \(L^* E\) to one that acts by isometries.

In the appendix we shall consider this alteration of the structure
group in more detail.

\subsection{Properties of the Polynomial Bundle}
\label{sec:polprop}

The construction of the polynomial loop bundle started from a
connection on the original bundle over \(M\).  However, it only
actually used the map \(\psi_X^0 : X^0 \to P X\) defined by the
parallel transport operator.  Thus as far as the polynomial loop
bundle is concerned, having a connection is overkill.  The connection
is useful, though, as it implies that the polynomial loop bundle came
from structure on the original manifold \(M\) and thus one can hope
for more structure on the polynomial loop bundle than has yet been
described.  In this section, we shall investigate this.  In the next,
we shall give an interpretation of the maps \(\psi_X^0\) in terms of
classifying maps and twisted K-theory.

Before examining the interesting properties of the polynomial loop
bundle, we list some basic ones that are fairly obvious:
\begin{proposition}
  Let \(M\) be a finite dimensional smooth manifold, \(E_1, E_2 \to
  M\) finite dimensional vector bundles over the same field with inner
  products and connections compatible with the inner products.
  \begin{enumerate}
  \item Let \(E = E_1 \oplus E_2\) orthogonally and equip \(E\) with
    the direct sum connection.  Then \(L_\pol E = L_\pol E_1 \oplus
    L_\pol E_2\).

  \item Suppose that \(E_1\) is real, then \(L_\pol (E_1 \otimes \C) =
    (L_\pol E_1) \otimes \C\).
    
  \item Suppose that \(E_1\) is complex, then \(L_\pol ({E_1}_\R) =
    (L_\pol E_1)_\R\).

  \item Let \(\psi : E_1 \to E_2\) be a bundle isomorphism which
    preserves the inner products and connections.  Then \(\psi\)
    defines an isomorphism \(L_\pol \psi : L_\pol E_1 \to L_\pol
    E_2\).

  \item Suppose that \(E_1\) with its inner product is a sub-bundle of
    \(E_2\) and that the covariant differential operator on \(E_1\) is
    of the form \(p \nabla\) where \(p : E_2 \to E_1\) is the
    orthogonal projection and \(\nabla\) is the covariant differential
    operator on \(E_2\).  Then it is not necessarily the case that
    \(L_\pol E_1 = L_\pol E_2 \cap L E_1\).
  \end{enumerate}
\end{proposition}

\begin{proof}
  Only the last of these is not immediate from the construction.  Let
  \(E_2\) be the bundle \(S^1 \times \C^2\) and \(E_1\) the bundle
  \(S^1 \times \C^1\).    Include \(E_1\) in \(E_2\) via the map \((t,
  1) \to (t, \frac1{\sqrt{2}} ( 1, e^{2 \pi i t}) )\).

  The loop space of \(E_1\) is \(L S^1 \times L \C\) and of \(E_2\) is
  \(L S^1 \times L \C^2\).  The polynomial loop space of \(E_2\) is
  \(L S^1 \times L_\pol \C^2\).  The inclusion \(L E_1 \to L E_2\) is
  given by:
  \[
  (\gamma, \beta) \to (\gamma, \frac1{\sqrt{2}} (\beta, e^{2 \pi i
  \gamma(t)} \beta)).
  \]
  Therefore \(L E_1 \cap L_\pol E_2\) consists of those loops
  \(\beta\) such that both \(\beta\) and \(e^{2 \pi i \gamma} \beta\)
  are polynomials.  We can choose \(\gamma\) such that whenever
  \(\beta\) is polynomial then \(e^{2 \pi i \gamma} \beta\) is not.
  Hence there is some \(\gamma\) such that above \(\gamma\) the
  fibres of \(L E_1\) and \(L_\pol E_2\) intersect trivially.
\end{proof}

The advantage of having the polynomial structure defined using a
connection on the original bundle is the relationship with the
diffeomorphism group of the circle.  For \(\sigma : S^1 \to S^1\)
smooth (not necessarily a diffeomorphism), \(\gamma : S^1 \to M\), and
\(\alpha \in L_\gamma E\), the following is a simple application of
the chain rule:
\begin{equation}
\label{eq:condiff}
D_{\gamma \circ \sigma} (\alpha \circ \sigma) =  \left((D_\gamma
\alpha) \circ \sigma \right) \sigma',
\end{equation}
where \(\sigma' : S^1 \to \R\) is such that \(d \sigma (\diff{}{t}) =
\sigma' \diff{}{t}\).

From this formula, two results can be derived:
\begin{proposition}
  \begin{enumerate}
  \item The action of \(\Diff(S^1)\) on \(L E\) does not preserve the
    sub-bundle \(L_\pol E\).  The subgroup of \(\Diff(S^1)\) which does
    preserve the sub-bundle \(L_\pol E\) is \(S^1 \rtimes \Z/2\) where
    the non-trivial element in the \(\Z/2\)-factor is the
    diffeomorphism \(t \to -t\).

  \item Let \(\nabla^a\) and \(\nabla^b\) be two different connections
    on \(E\).  The two polynomial bundles so defined are different.
  \end{enumerate}
\end{proposition}

\begin{proof}
  We shall consider the complex case so that we may talk about
  eigenvectors and eigenvalues of \(D_\gamma\).  The real case may be
  deduced from this.

  \begin{enumerate}
  \item For this, consider the situation over a constant loop.  There,
      \(L E\), resp.\ \(L_\pol E\), is \(E \otimes L \C\), resp.\ \(E
      \otimes L_\pol \C\).  The action of \(\Diff(S^1)\) on \(L E\) is
      given by its action on \(L \C\).  Thus if \(\sigma \in
      \Diff(S^1)\) preserves \(L_\pol E\) then it must preserve
      \(L_\pol \C\) within \(L \C\).

      The map \(t \to e^{2 \pi i t}\) lies in \(L_\pol \C\).  It is
      also the identification of \(S^1\) with \T.  Under \(\sigma\)
      this transforms to \(t \to e^{2 \pi i \sigma(t)}\).  As
      \(\sigma\) is a diffeomorphism of \(S^1\), this map must still
      be an identification of \(S^1\) with \T.  The only polynomials
      which do this are those of the form \(t \to \nu e^{\pm 2 \pi i
      t}\) for \(\nu \in \T\).  Hence if \(\sigma \in \Diff(S^1)\)
      preserves \(L_\pol E\) within \(L E\) then \(\sigma \in S^1
      \rtimes \Z/2\).

      The converse is direct from the equation~\ref{eq:condiff} since
      if \(\sigma \in S^1 \rtimes \Z/2\) then \(\sigma' = \pm 1\) so:
      \[
      D_{\gamma \circ \sigma} (\alpha \circ \sigma) = \pm
      (D_\gamma \alpha) \circ \sigma.
      \]
      Hence \(\sigma\) maps eigenvectors of \(D_\gamma\) to
      eigenvectors of \(D_{\gamma \circ \sigma}\) and thus preserves
      \(L_\pol E\).

    \item As \(\nabla^a\) and \(\nabla^b\) are different, there will
      be some loop \(\gamma\) such that \(D_\gamma^a\) and
      \(D_\gamma^b\) differ.  The difference will be a section
      \(\Phi\) of the bundle \(\mf{u}(\gamma^* E) \to S^1\), in other
      words an element of \(L_\gamma \mf{u}(E)\).

      If \(L_{\pol,\gamma}^a E = L_{\pol, \gamma}^b E\) then both are
      preserved under \(D_\gamma^a\) and \(D_\gamma^b\), hence under
      their difference.  Thus \(\Phi\) must be an element of \(L_\pol
      \mf{u}(E)\).

      By examining equation~\ref{eq:condiff}, we see that under the
      action of a smooth self-map \(\sigma\) of the circle, \(\Phi\)
      transforms to \((\Phi \circ \sigma) \sigma'\).  It is then a
      simple matter to find \(\sigma\) such that this is no longer a
      polynomial.  Hence even if we were unlucky enough initially to
      choose a loop \(\gamma\) with \(L_{\pol,\gamma}^a E =
      L_{\pol,\gamma}^b E\) then we can find some other loop \(\gamma
      \circ \sigma\) over which the fibres of the polynomial bundles
      differ. \qedhere
  \end{enumerate}
\end{proof}

It is straightforward to show that the result about the action of
\(\Diff(S^1)\) on \(L_\pol E\) generalises to the statement that the
subgroup of \(\Diff(S^1)\) which preserves \(L^? E\) is \(\Diff(S^1)
\cap L^? \C\) where the ``\(?\)'' represents some class of regularity
of loop.

In the light of this result, it is perhaps surprising that there is an
action of \(\Diff(S^1)\) on \(L_\pol E\) which covers the standard
action of \(\Diff(S^1)\) on \(L M\).  This comes about because the
\(\Diff(S^1)\)-action preserves the parallel transport operator.
Since all else was derived from that, we can make \(\Diff(S^1)\) act
on \(L_\pol E\).

We start with the group \(\Diff_0^+(S^1)\) of orientation and
basepoint preserving diffeomorphisms.  Since the whole diffeomorphism
group is the semi-direct product of this with \(S^1 \rtimes \Z/2\), an
action of this group together with the above action of \(S^1 \rtimes
\Z/2\) will give an action of the whole diffeomorphism group.

An element of \(\Diff_0^+(S^1)\) lifts canonically to an element of
\(\Diff_0^+(\R)\).  The image consists of those diffeomorphisms of \R
which satisfy \(\sigma(t + 1) = \sigma(t) + 1\).  This allows
\(\Diff_0^+(S^1)\) to act on paths as well as loops.

Let \(\sigma \in \Diff_0^+(S^1)\).  Recall that the bundle \(P^{M,0} E
\to L M\) has fibre \(P^{M,0}_\gamma E = P( E_{\gamma(0)})\).  Thus as
\(\gamma \circ \sigma(0) = \gamma(0)\), the bundles \(P^{M,0} E\) and
\(\sigma^*(P^{M,0} E)\) are genuinely the same bundle.  The bundle
\(P^L E\), meanwhile, has fibre \(P^L_\gamma E = \Gamma(\gamma^* E)\).
Thus there is a natural isomorphism \(P^L E \to \sigma^*(P^L E)\)
given by \(\alpha \to \alpha \circ \sigma\).

With these two isomorphisms, the square:
\[
\begin{CD}
  P^{M,0} E @>\Psi_E>> P^L E \\
  @| @V\sigma VV \\
  P^{M,0} E @>\Psi_E>> P^L E
\end{CD}
\]
does not commute.  To make it commute, we  need to transfer one action
of \(\sigma\) from one side to the other.  Clearly, the action of
\(\sigma\) on \(P^L E\) restricts to the standard action on \(L E\)
which we know does not preserve \(L_\pol E\).

It is also true that the action of \(\sigma\) on \(P^{M,0} E\)
preserves \(L^{M,0,h_E^{-1}} E\) and \(L^{M,0,h_E^{-1}}\).  Thus is
because the holonomy operator \(h_E\) is equivariant under the action
of \(\Diff_0^+(S^1)\).  Therefore, the action of \(\sigma\) on
\(P^{M,0} E\) when transferred to \(P^L E\) also restricts to an
action on \(L E\) and on \(L_\pol E\).

In formul\ae, the two actions of \(\Diff^+_0(S^1)\) are as follows:
any element of \(P_\gamma E\) can be written as \(\sum_j f^j \psi^0_E
v_j\) where \(\{v_1, \dotsc, v_n\}\) is a basis for
\(E_{\gamma(0)}\).  The usual action is:
\[
\sigma \left( \sum_j f^j \psi^0_E v_j \right) = \sum_j f^j \circ
\sigma \psi^0_E v_j
\]
and the new action is:
\[
\sigma \left( \sum_j f^j \psi^0_E v_j \right) = \sum_j f^j \psi^0_E
v_j.
\]

One way to make the distinction between the two actions is to have two
views of the bundle \(L E \to L M\).  In one, a fibre \(L_\gamma E\)
is inextricably linked to the points of \(\gamma(S^1)\).  In the
other, the fibre \(L_\gamma E\) is linked only to the map \(\gamma\).
In the former, reparametrising the loop \(\gamma\) does not change
\(\gamma(S^1)\) and so the fibres \(L_\gamma E\) and \(L_{\gamma \circ
\sigma} E\) are closely related.  Any reasonable -- in this view --
group action must preserve this relationship.  In the latter view,
reparametrising the loop \(\gamma\) changes it and so there is no
intrinsic relationship between the fibres \(L_\gamma E\) and
\(L_{\gamma \circ \sigma} E\).  Therefore there is no special
relationship for a reasonable group action to preserve.

\subsection{Loop Bundles and Twisted K-Theory}
\label{sec:twist}

As mentioned in the previous section, the construction of the
polynomial loop bundle started from a connection on the original
bundle over \(M\) but a connection provides rather more structure than
is needed.  The vital piece was the section of the bundle \(h_E\) of
\(Q^{\ad,0} \to L M\) and the isomorphism of \(L E\) with
\(L_\per^{M,0,h_E^{-1}} E\).  Thus any pair \((\chi, \Psi)\) where
\(\chi\) is a section of \(Q^{\ad,0} \to L M\) and \(\Psi\) is a
bundle isomorphism of \(L_\per^{M,0,\chi} E\) with \(L E\) will do.
However, if one wants \(\chi\) and \(\Psi\) to come from structure on
\(M\), a connection is the simplest starting point.  This ensures that
the fibres of the polynomial loop bundle are related to the points on
\(M\) over which they lie.

The section \(h_E\) (rather, \(h_E^{-1}\)) can be thought of as a type
of classifying map of the bundle \(L E\).  It is not, strictly
speaking, a classifying map as it does not land in \(B L G\).  Rather
it classifies \(L E\) ``up to constant loops''.  The basic idea of
this viewpoint is that when considering infinite dimensional geometry,
anything finite dimensional is relatively uninteresting or already
well-understood.  Therefore, saying that a bundle is trivial ``up to
constant loops'' is saying that it is really a finite dimensional
object that has been enhanced in some trivial way to make it appear
infinite dimensional and therefore is of little interest.  For
example, with polynomial loops, the ``polynomial'' part is really
defined for based loops.  To get the free polynomial loops, one simply
includes the constant loops in an appropriate way.

To make this slightly more mathematical, recall that the free loop
group, \(L G\), is the semi-direct product of the based loop group and
the constant loops.  That is, there is a split short exact sequence:
\[
\Omega G \to L G \leftrightarrow G.
\]
Thus within the class of \(L G\)-objects are those which come from
\(G\)-objects via the inclusion \(G \to L G\).  For example, within
the class of vector bundles over a space \(Y\) with fibre \(L \C^n\)
lie the vector bundles of the form \(E \otimes L \C\) for some
\(n\)-dimensional vector bundle \(E \to Y\).

There is a similar sequence of classifying spaces.  A particular
choice of classifying spaces, used for example in~\cite{rcas}, is:
\[
G \to E G \times_{\text{conj}} G \to B G.
\]

Given a classifying map \(Y \to B G\) we can thus pull-back the
\(G\)-bundle (which is not a principal bundle but rather a bundle of
groups) over \(Y\).  We interpret this as a bundle over \(Y\) with
fibre \(B \Omega G\).  Thus a section of this bundle defines a
\emph{twisted} principal \(\Omega G\)-bundle over \(Y\).  A section of
the \(B \Omega G\)-bundle is also a map from \(Y\) to \(E G
\times_{\text{conj}} G = B L G\) and thus classifies a principal \(L
G\)-bundle.

Conversely, given a classifying map \(Y \to B L G = E G
\times_{\text{conj}} G\), we can project down onto \(B G\) and
pull-back the \(G\)-bundle as above.  The original classifying map
then defines a section of this \(G\)-bundle and so a twisted principal
\(\Omega G\)-bundle over \(Y\).

Hence a principal \(L G\)-bundle can be interpreted as a principal
\(\Omega G\)-bundle twisted by a principal \(G\)-bundle.  The
\(G\)-bundle that defines the twisting is the pull-back of the
\emph{principal} \(G\)-bundle from \(B G\).  The \(B \Omega G\)-bundle
used above is the adjoint bundle of this principal \(G\)-bundle.

This is actually an unstable phenomenon, at least in the case of
\(U_n\).  There is a group homomorphism \(L U_n \to U_n \times
\gr_\res(H)\) which is a homotopy equivalence in the stable range.
Thus \(E U_\infty \times_{\text{conj}} U_\infty \simeq B U_\infty
\times U_\infty\).  In fact, by choosing appropriate models for
\(U_\infty\) and \(B U_\infty\) we can make this a homeomorphism.  Let
\(H\) be a complex, separable infinite dimensional Hilbert space.  Let
\(U(H)\) be the unitary operators on \(H\).  Let \(U_\m{K}\) be the
unitary operators on \(H\) of the form \(1 + T\) for some compact
operator \(T\).  Then by~\cite{nk} and~\cite{rp}, \(U(H)\) is
contractible and \(U_\m{K} \simeq U_\infty\).  The maps \(U(H)
\times_{\text{conj}} U_\m{K} \leftrightarrow B U_\m{K} \times
U_\m{K}\) are:
\[
[q,p] \to ( [q], q p q^{-1} ), \quad ([q], p) \to [q, q^{-1} p q],
\]
where we use the fact that \(U_\m{K}\) is normal in \(U(H)\).  Thus
the twisting of an \(\Omega U_\infty\)-bundle by a \(U_\infty\)-bundle
is trivial.

We end by noting in passing that the
splitting of \(B L U_\infty\) is related to the fact that \(L
U_\infty\) is its own (based) loop space and hence its own classifying
space.  Thus, for example, it defines a ring spectrum and hence a
generalised cohomology theory of period \(1\).  It is not a very
interesting theory as it is just \(K^0 + K^{-1}\).

\newpage
\section{The Dirac Operator on the Loop Space}
\label{sec:diracloop}

In this section we construct the Dirac operator on the loop space of
an appropriate manifold.  We start with a review of the theory of spin
in infinite dimensions and its links to loop groups.  We then turn to
the question of what structure on the original manifold gives rise to
a spin structure on the loop space.  Finally, we construct the Dirac
operator.

\subsection{Spin Structures and Polarisations}
\label{sec:spinpol}

In this section we shall review the essential details of the
construction of the spin representation in infinite dimensions, also
referred to as the Fock representation.  This is gleaned mostly from
\cite{rppr} with the application to loop spaces coming from
\cite{apgs}.

Let \(V\) be an infinite dimensional real vector space with a
continuous inner product, \(\ipvc\).  Let \(J\) be a choice of unitary
structure on \(V\); that is, \(J\) is an orthogonal transformation on
\(V\) such that \(J^2 = -1\).  Let \(V_J\) denote \(V\) with this
complex structure and let \(\ipc\) be the hermitian inner product on
\(V_J\) defined by \(\ip{u}{v} = \ipv{u}{v} + i \ipv{u}{J v}\).

Let \(\bbh_J\) be the Hilbert space completion of \(\exterior^\bullet
V_J\), the exterior power of \(V_J\), with respect to the inner
product:
\[
\ip{u_1 \wedge \dotsb \wedge u_k}{v_1 \wedge \dotsb \wedge v_l} =
\begin{cases}
0 & l \ne k \\ \det ( \ip{u_i}{v_j} ) & l = k.
\end{cases}
\]
Recall that \(\m{L}(\bbh_J)\) is the Banach space of (complex)
continuous linear maps from \(\bbh_J\) to itself.  Define operators
\(c : V \to \m{L}(\bbh_J)\) and \(a : V \to \m{L}(\bbh_J)\) by:
\begin{align*}
c(v)u_1 \wedge \dotsb \wedge u_k &= v \wedge u_1 \wedge \dotsb \wedge
u_k \\
a(v)u_1 \wedge \dotsb \wedge u_k &= \sum_{j=1}^k (-1)^{j-1}
\ip{u_i}{v} u_1 \wedge \dotsb \wedge \widehat{u_i} \wedge \dotsb
\wedge u_k.
\end{align*}
Let \(\pi : V \to \m{L}(\bbh_J)\) be the operator \(c + a\).

\begin{proposition}
The operator \(c\) is complex linear and \(a\) is conjugate linear,
regarding \(V\) as \(V_J\), and they satisfy the \emph{canonical
anti-commutation relations}:
\begin{align*}
\left\{c(u),a(v)\right\} &= \ip{u}{v} \\
\left\{c(u),c(v)\right\} = \left\{a(u), a(v)\right\} &= 0
\end{align*}
where for operators \(X,Y\), \(\left\{X,Y\right\} = X Y + Y X\).

Hence \(\pi\) is real linear and satisfies \(\pi(v)^2 = \ipv{v}{v}
I\).
\end{proposition}

The map \(\pi\) is called \emph{Clifford multiplication}.  The space
\(\bbh_J\) decomposes as \(\bbh_J^+ \oplus \bbh_J^-\) with
\(\bbh_J^+\) the completion of \(\exterior^{\text{ev}} V_J\) and
\(\bbh_J^-\) of \(\exterior^{\text{odd}} V_J\).  With respect to this
grading, Clifford multiplication is of odd degree.  That is, it
interchanges the factors.

It is fairly obvious, and is described in~\cite[theorem 1.2.7]{rppr},
that the construction of the Fock representation factors through the
Hilbert completion of \(V\) defined by the inner product.  Thus if
\(\widetilde{V}\) is a subspace of \(V\), possibly with a finer
topology, that is dense in \(V\) with the inner product topology and
such that \(J\) restricts to a unitary structure on \(\widetilde{V}\)
then the Fock representations of \((V, J)\) and \((\widetilde{V}, J)\)
are the same.

The \emph{implementation question} is the following: let \(O(V)\) be
the orthogonal group of \(V\).  For which \(g \in O(V)\) is there some
\(U_g \in U(\bbh_J)\) such that \(\pi(g v) = U_g \pi(v) U_g^{-1}\)?
It is answered by:

\begin{theorem}[{\cite[ch 3]{rppr}}]
For \(g \in O(V)\) there is some \(U_g \in U(\bbh_J)\) such that
\(\pi(g v) = U_g \pi(v) U_g^{-1}\) if and only if \(\left[g,
J\right]\) is a Hilbert-Schmidt operator.  Moreover, if \(U_g\) and
\(U_g'\) both implement \(g\) then \(U_g = \lambda U_g'\) for some
\(\lambda \in S^1\).
\end{theorem}

An operator \(T : H_1 \to H_2\) between Hilbert spaces is said to be
\emph{Hilbert-Schmidt} if for some, and hence every, orthogonal basis
\(\{e_i\}\) of \(H_1\) then \((\norm[T e_i])\) is square summable.
The subgroup of \(O(V)\) consisting of \(g\) such that
\(\left[g,J\right]\) is Hilbert-Schmidt is written \(O_J(V)\) in
\cite{rppr}.

An intimately related problem is that of \emph{equivalence}: given
unitary structures \(J\) and \(K\) on \(V\), what condition is
equivalent to there being a unitary transformation \(T : \bbh_J \to
\bbh_K\) such that \(\pi_J(v) = T \pi_K(v) T^{-1}\)?  The answer is
given by:

\begin{theorem}[{\cite[ch 3]{rppr}}]
Let \(J\) and \(K\) be unitary structures on \(V\).  The Fock
representations \(\bbh_J\) and \(\bbh_K\) are unitarily equivalent if
and only if \(J - K\) is Hilbert-Schmidt.
\end{theorem}

Thus one could define a \emph{Fock structure} on \(V\) to be an
equivalence class of unitary structures.  The Fock representation
would then only depend on this class, rather than the explicit choice
of unitary structure.  This idea provides a neat sequ\'e{} from the
theory of Fock representations to that of polarisations.

There are various equivalent definitions of a polarisation, we choose
the one that is closest to the theory of unitary structures.  The
theory of polarisations and the relationship with loop groups is the
subject of \cite{apgs}.  The following definitions are equivalent to
those from \cite[ch 6]{apgs} although we have used notation similar to
that of \cite{rppr} for better comparison with the theory of unitary
structures.

\begin{defn}
Let \(H\) be a complex Hilbert space.  A \emph{polarising operator} on
\(H\) is an operator \(J \in \m{L}(H)\) such that \(J^2 + I\) is trace
class and \(J \pm i I\) are not finite rank.

A \emph{polarisation} on \(H\) is an equivalence class of polarising
operators under the relation \(J_1 \sim J_2\) if and only if \(J_1 -
J_2\) is Hilbert-Schmidt.

Let \(\m{J}\) be a polarisation on \(H\).  The \emph{restricted
general linear group} of \(H\) with respect to \(\m{J}\),
\(\gl_\m{J}(H)\), is defined as the subgroup of \(\gl(H)\) consisting
of those \(A\) for which \(\left[A,J\right]\) is Hilbert-Schmidt for
one, and hence all, \(J \in \m{J}\).
\end{defn}

In \cite{apgs}, the notation used is \(\gl_\res(H)\).  The notation
\(\gl_\m{J}(H)\) emphasises the dependence on the polarisation
\(\m{J}\).  The operator used in the above definition is slightly
different from the operator \(J\) used in \cite[ch 6]{apgs}.  To get
from the one to the other, multiply by \(-i\).

Clearly a polarising operator \(J\) defines a polarisation by taking
the equivalence class of \(J\).  Thus a unitary structure \(J\) on a
real Hilbert space \(H\) gives rise to a polarisation on the
complexification \(H_\C\) by taking the equivalence class of \(J\),
extended to the complexification by linearity.  With respect to this
polarisation, it is evident that \(O_J(H) = \gl_\m{J}(H_\C) \cap
O(H)\).

There are three equivalent definitions of a unitary structure given in
\cite[ch 2.1]{rppr}.  Using these correspondences, a careful
examination of \cite[ch 12]{apgs} reveals that the standard unitary
structure on \(L^2(S^1,\R^{2n})\) is defined in the following way: Let
\(\{e_k\}\) be the standard basis for \(\R^{2n}\).  Let \(J_0 :
\R^{2n} \to \R^{2n}\) be the complex structure \(J_0 e_{2k} = e_{2k -
1}\), \(J_0 e_{2k - 1} = - e_{2k}\).  The unitary structure on
\(L^2(S^1, \R^{2n})\) is defined by the operator \(J\) which
satisfies:
\begin{align*}
J( v \cos k\theta ) &= v \sin k\theta \\ J(v \sin k\theta ) &= -v \cos
k\theta \\ J(v) &= J_0(v).
\end{align*}
Here we identify \(\R^{2n}\) with the subspace of constant loops in
\(L^2(S^1, \R^{2n})\).

The standard polarisation operator \(J\) on \(L^2(S^1, \C^m)\)
satisfies the identity:
\[
J(v z^k) = - (-1)^{\text{sign}(k)} i v z^k.
\]

\begin{proposition}
  \label{prop:unipol}
  The standard polarisation on \(L^2(S^1, \C^{2n})\) is that defined
  by the standard unitary structure on \(L^2(S^1, \R^{2n})\).  If
  \(m\) is odd, the standard polarisation on \(L^2(S^1, \C^m)\) does
  not contain a unitary structure for \(L^2(S^1, \R^m)\).
\end{proposition}

\begin{proof}
  To distinguish the operators, let \(J_\R\) denote the unitary
  structure on \(L^2(S^1,\R^{2n})\) and also its extension to
  \(L^2(S^1,\C^{2n})\).  Let \(J_\C\) be the polarising operator on
  \(L^2(S^1,\C^m)\).  The first part of the proposition follows from
  the observation that \(J_\R\) and \(J_\C\) agree on the subspace of
  \(L^2(S^1,\C^{2n})\) consisting of loops orthogonal to the constant
  loops.  This has finite codimension and so \(J_\C - J_\R\) is finite
  rank.  Thus \(J_\R\) and \(J_\C\) define the same polarisation on
  \(L^2(S^1,\C^{2n})\).

  Let \(m\) be odd.  Let \(\{e_k\}\) be the standard basis for
  \(\R^m\).  Let \(J_0 : \R^m \to \R^m\) be the map \(J_0 e_{2k} =
  e_{2k-1}\), \(J_0 e_{2k-1} = -e_{2k}\), \(J_0 e_m = 0\).  Let
  \(J_\R\) be the map on \(L^2(S^1, \R^m)\) defined using \(J_0\) as
  for the even dimensional case.  This restricts to a unitary
  structure on the subspace \(\left\langle e_m \right \rangle^\perp\).
  As before, \(J_\R\) and \(J_\C\) agree on the subspace of loops
  orthogonal to the constant loops and thus define the same
  polarisation on \(L^2(S^1,\C^m)\).

  Let \(K\) be a unitary structure on \(L^2(S^1, \R^m)\).  The space
  \(L^2(S^1, \C^m)\) decomposes orthogonally according to the
  eigenspaces of \(J_\R\) and of \(K\).  Corresponding to \(J_\R\) we
  have \(L^2(S^1, \C^m) = V_+ \oplus V_- \oplus \C\) as \(\pm
  i\)-eigenspaces and the \(0\)-eigenspace.  Corresponding to \(K\) we
  have \(L^2(S^1, \C^m) = W_+ \oplus W_-\).  Let \(\Sigma\) denote the
  operation of complex conjugation on \(L^2(S^1, \C^m)\).  Then
  \(\Sigma W_\pm = W_\mp\), \(\Sigma V_\pm = V_\mp\), and \(\Sigma \C
  = \C\).

  The identity map decomposes as the matrix:
  \[
  \begin{bmatrix}
    a & b & c \\ d & e & f
  \end{bmatrix}
  : V_+ \oplus V_- \oplus \C \to W_+ \oplus W_-.
  \]
  Here \(a : V_+ \to W_+\) is the inclusion of \(V_+\) followed by the
  projection onto \(W_+\), and similarly for the other entries.  Since
  the identity map commutes with complex conjugation, \(d = \Sigma b
  \Sigma\), \(e = \Sigma a \Sigma\), and \(f = \Sigma c \Sigma\).

  Now assume that \(J_\R - K\) is Hilbert-Schmidt.

  The operator \(b : V_- \to W_+\) can be written as \(\frac14(I - i
  K)(I + i J_\R)P\) where \(P : V_+ \oplus V_- \oplus \C \to V_+
  \oplus V_-\) is the orthogonal projection.  This expands to
  \(\frac14(I + K J_\R + i(J_\R - K))P\).  As \(K^2 = -I\), \(I + K
  J_\R = K(J_\R - K)\) and therefore \(b\) is Hilbert-Schmidt.
  Similarly, \(d\) is Hilbert-Schmidt.  Since \(c\) and \(f\) have
  domain \C, they are finite rank.  Thus the operator \(a + e\)
  differs from the identity by a compact operator so is Fredholm of
  index zero.

  Since \(a\) and \(e\) start from orthogonal subspaces and end in
  orthogonal subspaces, the fact that \(a + e\) is Fredholm implies
  that both \(a\) and \(e\) are also Fredholm.  The identity \(e =
  \Sigma a \Sigma\) then implies that \(\ind a = \ind e\).  The matrix
  form of \(a + e\) is:
  \[
  \begin{bmatrix}
    a & 0 & 0 \\ 0 & e & 0
  \end{bmatrix}
  \]
  from which it is evident that the index of \(a + e\) is \(\ind a +
  \ind e + 1\).  This is incompatible with \(\ind a = \ind e\) and so
  we deduce that \(J_\R - K\) cannot be Hilbert-Schmidt.  Hence there
  is no unitary structure for \(L^2(S^1, \R^m)\) in the standard
  polarisation of \(L^2(S^1, \C^m)\).
\end{proof}

For the record, we note the following properties of the groups
associated to the standard polarisation on \(L^2(S^1, \C^{2n})\) and
the standard unitary structure on \(L^2(S^1, \R^{2n})\).

\begin{lemma}
Let \(H = L^2(S^1,\R^{2n})\) and let \(J\) be the standard unitary
structure on \(H\).  Let \(H_\C = L^2(S^1,\C^{2n})\) be the
complexification and \(\m{J}\) the standard polarisation on \(H_\C\).

\begin{enumerate}
\item \(O_J(H) = \gl_\m{J}(H_\C) \cap O(H)\);
\item let \(U_\m{J}(H_\C) = \gl_\m{J}(H_\C) \cap U(H_\C)\), then
  \(U_\m{J}(H_\C) \to \gl_\m{J}(H_\C)\) is a deformation retract;
\item let \(\gl_J(H) = \gl_\m{J}(H_\C) \cap \gl(H)\), then \(O_J(H)
  \to \gl_J(H)\) is a deformation retract; and
\item \(U_\m{J}(H_\C) \simeq \Omega U\), \(O_J(H) \simeq \Omega O\).
\end{enumerate}
\end{lemma}

In \cite[ch 6]{apgs}, it is shown that the natural action of
\(LU_{2n}\) on \(H_\C\) defines an inclusion \(LU_{2n} \to
U_\m{J}(H_\C)\).  Since \(LO_{2n} = LU_{2n} \cap O(H)\), it follows
that the natural action of \(LO_{2n}\) on \(H\) defines an inclusion
\(LO_{2n} \to O_J(H)\).

The action of \(O_J(H)\) on \(\bbh_J\) is projective.  That is, there
is a central \(S^1\)-extension of \(O_J(H)\), usually written
\(\pin_J(H)\) (the identity component being \(\spin_J(H)\)), which
acts unitarily on \(\bbh_J\).  This central extension is classified by
a generator of \(H^2(O_J(H),\Z)\), which is isomorphic to \Z.

Examining \(LO_{2n}\), we see that it has four components.  The
identity component is the semi-direct product \(SO_{2n} \times \Omega
\spin_{2n}\) which has double cover \(L\spin_{2n}\).  The central
extension of \(O_J(H)\) pulls back to a central \(S^1\)-extension of
\(L\spin_{2n}\) written \(\widetilde{L}\spin_{2n}\).  This is
classified by a generator of \(H^2(L\spin_{2n}, \Z)\), which is also
isomorphic to \Z.  Note also that the transgression map \(\tau :
H^\bullet(\spin_{2n}, \Z) \to H^{\bullet - 1}(L\spin_{2n}, \Z)\) is an
isomorphism from degree \(3\) to degree \(2\).

The two components of \(\pin_J\) can be easily distinguished.  Recall
that \(\bbh_J\) decomposes as \(\bbh_J^+ \oplus \bbh_J^-\).  The
identity component of \(\pin_J(H)\), whence also
\(\widetilde{L}\spin_{2n}\), preserves this decomposition.  The other
component swaps the factors.

Finally, the circle action on \(L^2(S^1,\R^{2n})\) lies in \(O_J(H)\)
and has a canonical lift to \(\pin_J(H)\).  This defines a circle
action on \(\bbh_J\).  The circle action on \(L\spin_{2n}\) therefore
lifts to \(\widetilde{L}\spin_{2n}\) and the action of
\(\widetilde{L}\spin_{2n}\) on \(\bbh_J\) is circle equivariant.

\subsection{String Manifolds and Spin Connections}
\label{sec:string}

In this section we explain how a string structure on a manifold
defines a connection on the spin bundle of the loop space.  Let \(M\)
be an oriented, Riemannian manifold of even dimension \(d\).  Let \(R
\to M\) be the principal \(SO_d\)-bundle determined by the metric and
the orientation.  Let \(\omega : T R \to \mf{so}_d\) be the
Levi-Civita connection on \(M\).

The group \(\spin_d\) is the connected double cover of \(SO_d\),
universal if \(d > 2\).  A \emph{spin structure} on \(M\) is a
principal \(\spin_d\)-bundle \(Q \to M\) such that \(Q\) is a double
covering of \(R\) and the following diagram commutes:
\[
\begin{CD}
\spin_d \times Q @>>> Q \\ @VVV @VVV \\ SO_d \times R @>>> R.
\end{CD}
\]

The manifold \(M\) admits a spin structure if and only if \(w_2(M) =
0\); the set of isomorphism classes of spin structures is in bijective
correspondence with \(H^1(M; \Z_2)\).

In order that the loop space, \(L M\), admit a spin structure the
structure group of \(L M\) must lift from \(L\spin_d\) to
\(\widetilde{L} \spin_d\).  We would also like this to be
\(S^1\)-equivariant.  The \(L\spin_d\)-principal bundle on \(L M\) is
\(L Q\).  Thus we are asking for an \(S^1\)-bundle, equivalently a
line bundle, over \(L Q\) with certain properties.  The primary
property is that on fibres it must pull-back to the fibration \(S^1
\to \widetilde{L} \spin_d \to L \spin_d\).

As explained in \cite[ch VI]{jb}, line bundles on loop spaces are
closely related to gerbes on the original manifold.  In particular,
the central extension \(\widetilde{L} \spin_d\) of \(L\spin_d\)
corresponds to the gerbe of \(\spin_d\) classified by the generator of
\(H^3(\spin_d; \Z)\) (recall that as a simply connected, simple Lie
group, there is a canonical isomorphism of \(H^3(\spin_d; \Z)\) with
\Z and hence a canonical generator).  Rather than asking for a line
bundle over \(L Q\) we therefore ask for a gerbe over \(Q\).  This has
the considerable advantage that the line bundle defined by the gerbe
will be \(\Diff^+(S^1)\)-equivariant.

We have to answer the following question: what is the obstruction to
constructing a gerbe on \(Q\) which on fibres pulls-back to the
fundamental gerbe on \(\spin_d\)?  We can rephrase this question in
cohomological terms where it becomes: when can we find an element \(a
\in H^3(Q;\Z)\) such that if \(i : \spin_d \to Q\) is the inclusion of
a fibre then \(i^* a\) is the generator of \(H^3(\spin_d;\Z)\)?

To answer this we examine the Serre spectral sequence of the fibration
\(\spin_d \to Q \to M\).  The first part of the \(E_2\)-term is:

\begin{picture}(0,0)
\put(25,-75){\line(0,1){65}} \put(6,-62){\line(1,0){305}}
\end{picture}
\[
\begin{matrix}
3 && H^3(\spin_d; \Z) \\ 2 && 0 & 0 & 0 \\ 1 && 0 & 0 & 0 & 0 \\ 0 &&
H^0(M;\Z) & H^1(M;\Z) & H^2(M;\Z) & H^3(M;\Z) & H^4(M;\Z) \\[6pt] && 0
& 1 & 2 & 3 & 4
\end{matrix}
\]

This contains all the possible contributions to \(H^3(Q;\Z)\).  The
only part that might not persist to the \(E_\infty\)-term is
\(H^3(\spin_d;\Z)\) in the \((0,3)\) position.  This persists until
the \(E_4\)-term where the differential is \(d_4 : H^3(\spin_d;\Z) \to
H^4(M;\Z)\).  Let \(\lambda \in H^4(M;\Z)\) denote the image of the
canonical generator of \(H^3(\spin_d; \Z)\) under \(d_4\).  If
\(\lambda = 0\) then \(H^3(Q;\Z) \cong H^3(M;\Z) \oplus H^3(\spin_d;
\Z)\) and the inclusion of a fibre induces the projection \(H^3(M;\Z)
\oplus H^3(\spin_d; \Z) \to H^3(\spin_d; \Z)\).  If \(\lambda \ne 0\)
then \(H^3(Q; \Z) = H^3(M;\Z)\) and the inclusion of a fibre is the
zero map on \(H^3\).  The class \(\lambda\) is known to satisfy \(2
\lambda = p_1(M)\) which has led to it being written as \(p_1(M)/2\).
This notation is somewhat misleading as \(\lambda\) depends on the
choice of spin structure on \(M\).

\begin{defn}
A manifold \(M\) is a \emph{string} manifold if it is an oriented,
Riemannian, spin manifold such that \(\lambda = 0\) together with a
choice of \emph{string structure}.  That is, a choice of gerbe,
\(\m{G}\), over the spin structure \(Q \to M\) which on fibres is the
fundamental gerbe on \(\spin_d\).
\end{defn}

Once we have a string structure, there is a natural notion of a string
connection.

\begin{defn}
A \emph{string connection} on a string manifold with string manifold
with string structure \(\m{G}\) consists of the Levi-Civita connection
on \(Q\) and a \(\spin_d\)-equivariant connective structure on the
gerbe \(\m{G}\).
\end{defn}

\begin{theorem}
A string connection on \(M\) defines a \(\Diff^+(S^1)\)-equivariant
spin connection on \(L M\).
\end{theorem}

Compare this result with that of \cite{pm2}.

\begin{proof}
The Levi-Civita connection on \(M\) is a map \(\omega : T R \to
\mf{so}_d\).  As \(\spin_d \to SO_d\) is a covering map, it is a local
diffeomorphism and so \(\mf{spin}_d = \mf{so}_d\).  Thus the
Levi-Civita connection lifts to a connection on \(Q\) via \(\omega' :
T Q \to T R \xrightarrow{\omega} \mf{so}_d = \mf{spin}_d\).  The loop
of this is a \(\Diff^+(S^1)\)-equivariant map \(L\omega' : T L Q \to
L\mf{spin}_d\).  This is also a connection.

The gerbe with its connective structure defines a
\(\Diff^+(S^1)\)-equivariant \(S^1\)-bundle \(\widetilde{L} Q \to L
Q\) with a connection \(\alpha : T \widetilde{L} Q \to \R\).  As the
gerbe on \(M\) pulls back to the fundamental gerbe on fibres, so also
\(\widetilde{L} Q \to L Q\) pulls back to \(\widetilde{L}\spin_d \to L
\spin_d\) on fibres.  Also, the connection is \(L
\spin_d\)-equivariant.  Hence \(L \omega' \oplus \alpha : T
\widetilde{L} Q \to L\mf{spin}_d \oplus \R\) is a connection on
\(\widetilde{L} Q\).
\end{proof}

\subsection{The Dirac Operator}
\label{sec:dirac}

We can now construct the Dirac operator.  Let \(M\) be a finite
dimensional, simply connected, string manifold.  The loop space \(L
M\) thus has a spin structure with spin connection.  The Levi-Civita
connection on the tangent bundle of \(M\) defines the polynomial loop
bundle, \(L_\pol T M\), and thus the Hilbert completion of \(T^* L M =
L^* T M\).

The spinor bundles, \(S^+, S^-\), of \(L M\) are constructed from
\(L^2 T M\).  As \(L^2 T M\) is a real Hilbert bundle, it is
canonically isomorphic to its dual, \(L^{2,*} T M\).  Thus we can
view the spinor bundle as being constructed from either \(L^2 T M\) or
\(L^{2,*} T M\) as seems appropriate.  The point of having the two
approaches is to distinguish between \(L^2 T M\) as the completion of
\(L T M\) and of \(L^* T M\) and thus determine which of the finite
dimensional constructions we are generalising.

The spin connection on \(L M\) defines a covariant differential
operator:
\[
\nabla : \Gamma(S^+) \to \Gamma(\m{L}(T L M, S^+)).
\]

By taking the spinor bundles from \(L^{2,*} T M\), we are considering
this to be the completion of \(L^* T M = T^* L M\).  Thus we consider
Clifford multiplication to be a fibrewise map \(L^* T M \to
\m{L}(S)\).  The following proposition is essentially the
\emph{remarkable isomorphism} and will enable us to compose this with
the covariant differential operator to define the Dirac operator.

\begin{proposition}
\label{prop:cmext}
Let \(V\) be a complete nuclear reflexive space with a continuous
inner product.  Let \(J\) be a unitary structure on \(V\).  The map
\(\pi : V \to \m{L}(\bbh_J)\) defines a continuous linear map \(\pi :
\m{L}(V^*, \bbh_J) \to \bbh_J\).
\end{proposition}

Before proving this, we show how this leads to the definition of the
Dirac operator.  We are considering \(S^+\) and \(S^-\) to be
constructed from the cotangent bundle, \(T^* L M\).  Therefore, we
take \(V\) in the statement of the proposition to be \(L^* \R^n\)
which is a complete nuclear reflexive space.  Since its dual is \(L
\R^n\), Clifford multiplication defines a fibrewise linear map:
\[
\pi : \m{L}(L T M, S^+) \to S^-.
\]

\begin{defn}
  The Dirac operator, \(\dirac : \Gamma(S^+) \to \Gamma(S^-)\), on \(L
  M\) is the composition:
  \[
  \dirac : \Gamma(S^+) \xrightarrow{\nabla} \Gamma(\m{L}(T L M, S^+))
  \xrightarrow{\pi} \Gamma(S^-).
  \]
\end{defn}

Since every piece of structure that went into its construction is
equivariant under rotations of the circle, the Dirac operator is
similarly equivariant.

We conclude with the proof of proposition~\ref{prop:cmext}:

\begin{proof}[Proof of proposition~\nmref{prop:cmext}]
Let \(H\) denote the Hilbert space completion of \(V\) with respect to
the inner product topology defined by the inner product on \(V\).
From \cite[ch 2.4]{rppr}, we know that \(\pi : V \to \m{L}(\bbh_J)\)
extends to an isometric inclusion \(\pi : H \to \m{L}(\bbh_J)\).  The
map \(H \times \bbh_J \to \bbh_J\), \((x,\xi) \to \pi(x)\xi\), is
therefore continuous.  From \cite[ch III, \S6]{hs}, it extends to a
continuous linear map with domain the projective tensor product \(H
\wotimes \bbh_J\).

The inclusion \(V \to H\) induces a continuous linear map \(V
\wotimes \bbh_J \to H \wotimes \bbh_J\).  From
\cite[ch IV, \S 9.4]{hs}, as \(V\) is a complete nuclear space then
the space \(V \wotimes \bbh_J\) is isomorphic to
\(\m{L}_e(V_\tau^*, \bbh_J)\); where this denotes the space of linear
maps from \(V^*\) to \(\bbh_J\).  The topology on \(V^*\) is the
Mackay topology and the topology on the space of maps is that of
uniform convergence on equicontinuous sets.

From \cite[ch IV, \S 5]{hs} we deduce that as \(V\) is reflexive, the
Mackay topology on the dual agrees with the strong topology.  Also as
\(V\) is reflexive, it is barrelled and so equicontinuous sets in
\(V^*\) are the same as bounded sets.  Hence \(\m{L}_e(V_\tau^*,
\bbh_J) = \m{L}(V^*, \bbh_J)\).
\end{proof}

\newpage
\appendix

\addcontentsline{toc}{section}{Appendix: Inner Products on
  the Space of Distributions}
\section*{Appendix: Inner Products on the Space of Distributions}
\renewcommand{\thesection}{A}

In this appendix we examine inner products on \(L^*\R^n\).  The goal
is to classify the inner products on \(L^*\R^n\) which have the
following properties: the inner product is invariant under the circle
action, the involution of reversing loops is orthogonal, and the
operations of multiplication by \(\cos \theta\) and \(\sin \theta\)
are continuous.

We shall actually work with \(\m{S}^*\), the dual of the space of
rapidly decreasing, complex-valued, \Z-indexed sequences.  As a
sequence space, this is particularly simple to describe and therefore
to work with.  Taking Fourier coefficients defines an isomorphism
\(L\C \to \m{S}\) which allows us to transfer information from
\(\m{S}^*\) to \(L^*\C\).  Using the description of \(L\C^n\) as \(L\C
\otimes \C^n\), we can extend the description to the dual of
\(L\C^n\), and thence to the dual of the underlying real space
\(L\R^n\).

As a preliminary, we shall prove a negative result.  We shall show
that there is no ``natural'' inner product on \(L^*\C\).  That is, if
\(L\C^\times = L(\C^\times)\) denotes the space of never-zero smooth
loops in \C then there is no inner product on \(L^*\C\) such that the
group \(L\C^\times\) acts continuously with respect to the inner
product topology.  This is in stark contrast to the situation for
\(L\C\) where \(L\C^\times\) does act continuously with respect to the
standard inner product.

\begin{theorem}
Let \(L_c \C\) be a class of loops in \C with the following
properties:
\begin{enumerate}
\item there are continuous inclusions \(\C \to L_c \C \to L^{1,\infty}
  \C\), where \C corresponds to the constant loops and \(L^{1,\infty}
  \C\) is the space of continuously differentiable loops;

\item the class of loops is preserved under products; thus \(L_c
  \C^\times\) acts on \(L_c \C\) and hence, via the adjoint map, on
  \(L_c^* \C\);

\item \(L_c \C\) is reflexive;

\item \(L_c \C\) cannot be given the structure of a Hilbert space;
\end{enumerate}
then for any inner product on \(L_c^* \C\) there is some \(\alpha \in
L_c \C^\times\) which acts unboundedly on \(L_c^* \C\) with respect to
the inner product topology.
\end{theorem}

\begin{proof}
Let \(\ipc\) be a continuous inner product on \(L_c^* \C\).  Let \(H\)
denote the Hilbert space completion of \(L_c^* \C\) with respect to
\(\ipc\).  The dual of the inclusion of \(L_c^* \C\) in \(H\) is a map
\(H^* \to {L_c^*}^* \C = L_c \C\).
				   
Suppose that \(L_c \C^\times\) acts continuously on \(L_c^* \C\) with
respect to the inner product topology.  This implies that \(H^*\) is
preserved in \(L_c \C\) by \(L_c \C^\times\).  Suppose that \(H^* \cap
L_c \C^\times \ne \emptyset\).  Because \(L_c \C^\times\) is a group,
this implies that \(L_c \C^\times \subseteq H^*\).  The linear span of
\(L_c \C^\times\) is \(L_c \C\) so \(H^* = L_c \C\).  However, this
implies that \(H^* \to L_c \C\) is a continuous, linear bijection from
a Hilbert space onto \(L_c \C\) which contradicts the fourth
assumption.

Thus we need to show that the other assumptions imply that \(H^* \cap
L_c \C^\times \ne \emptyset\).  In other words, we need to show that
there is an element in \(H^*\) which is never zero.  To do this, we
shall use the Banach-Steinhaus theorem as stated in \cite[III, \S
4.6]{hs}.  As \(L_c \C\) is reflexive, it is the dual of \(L_c^* \C\).
We shall write the evaluation of \(\alpha \in L_c \C\) on \(a \in L_c^*
\C\) as \(a(\alpha)\) rather than \(\alpha(a)\) to avoid confusion
with the notation \(\alpha(\lambda)\) for the evaluation of \(\alpha\)
on \(\lambda \in S^1\).

From the corollary to \cite[IV, \S 2.3]{hs}, as the inclusion \(L_c^*
\C \to H\) is injective with weakly dense image, the map \(H^* \to L_c^*
\C\) is also injective with weakly dense image.  Thus there is a
sequence \((\alpha_n)\) in \(H^*\) which converges weakly to \(1\).
That is, for all \(a \in L_c^* \C\), \((a(\alpha_n))\) converges in \C
to \(a(1)\).  The space \(L_c^* \C\) is reflexive, hence barrelled, and
so the Banach-Steinhaus theorem applies.  This states that
\((\alpha_n)\) converges to \(1\) uniformly on each compact subset of
\(L_c^* \C\).  We shall find a particularly convenient compact subset
of \(L_c^* \C\).

The norm on \(L^{1,\infty} \C\) is \(\norm[\gamma]_{1,\infty} =
\sup\{\abs{\gamma(\lambda)}, \abs{\gamma'(\lambda)}\}\).  For
\(\lambda \in S^1\), there is an element \(e_\lambda\) of
\(L^{1,\infty} \C\) which evaluates a loop at time \(\lambda\).  If
\(\gamma \in L^{1,\infty} \C\) with \(\norm[\gamma]_{1,\infty} \le 1\)
then \(\gamma\) is Lipschitz with constant \(K \le 1\).  Therefore
\(\abs{e_\lambda(\gamma) - e_{\lambda'}(\gamma)}\) is less than or
equal to the smaller angle between \(\lambda\) and \(\lambda'\).
Hence \(\lambda \to e_\lambda\) is a continuous map from \(S^1\) to
\({L^{1,\infty}}^* \C\).  Composing this with the dual of the map \(L_c
\C \to L^{1,\infty} \C\) defines a continuous map \(S^1 \to L_c^*
\C\).  Its image is thus compact and therefore \((\alpha_n) \to 1\)
uniformly on \(\{e_\lambda : \lambda \in S^1\}\).

Hence there is some \(N\) such that for \(n \ge N\),
\(\abs{e_\lambda(\alpha_n) - e_\lambda(1)} < 1\) for all \(\lambda \in
S^1\).  Thus \(\abs{\alpha_N(\lambda) - 1} < 1\) so
\(\alpha_N(\lambda) \ne 0\) for all \(\lambda \in S^1\).  Hence \(H^*\)
contains an element which is never zero.
\end{proof}

\subsection{Inner Products on Distribution Space}

In this section we investigate those inner products on \(\m{S}^*\)
which, under the isomorphism \(\m{S}^* \cong L^* \C\), are invariant
under the circle action and the involution of reversing loops, and
such that multiplication by \(z\) is continuous in the inner product
topology.  In this investigation, we use \(\m{S}^*\) because it is a
sequence space and so we have a good presentation of elements of
\(\m{S}^*\) and of operators acting on it.  We start by transferring
the aforementioned operators from \(L^*\C\) to \(\m{S}^*\).

\begin{defn}
Define the operators \(R_\lambda\) for \(\lambda \in S^1\), \(\iota\),
and \(z\) on \(\m{S}\) to be the operators corresponding under the
Fourier isomorphism \(\m{S} \cong L\C\) to rotation by \(\lambda\),
reversal of the circle, and multiplication by \(z\), respectively.  We
shall use the same notation for their adjoints which act on
\(\m{S}^*\).
\end{defn}

The maps \(\lambda \to R_\lambda \in \m{L}(\m{S})\) and \(\lambda \to
R_\lambda \in \m{L}(\m{S}^*)\) define an action of the circle on
\(\m{S}\) and \(\m{S}^*\) respectively.  We shall refer to \(\iota\) as
the \emph{natural involution} on \(\m{S}\) and \(\m{S}^*\).

For \(p \in \Z\), let \(e^p \in \m{S}\) and \(e_p \in \m{S}^*\) both
denote the sequence with a \(1\) in the \(p\)th place and zero
elsewhere.  The sets \(\{e^p\}\) and \(\{e_p\}\) are topologically
free bases for \(\m{S}\) and \(\m{S}^*\) respectively.

\begin{lemma}
\label{lem:opprop}
In terms of the bases \(\{e^p\}\) and \(\{e_p\}\), the operators
\(R_\lambda\), \(\iota\), and \(z\) are given by the formul\ae:
\begin{align*}
R_\lambda e^p &= \lambda^p e^p &
\iota e^p &= e^{-p} &
z e^p &= e^{p+1} \\
R_\lambda e_p &= \lambda^{-p} e_p &
\iota e_p &= e_{-p} &
z e_p &= e_{p-1} 
\end{align*}
\end{lemma}

\begin{defn}
Let \(\m{C}\) denote the cone of positive semi-definite, sesquilinear
forms on \(\m{S}^*\) which are invariant under the action of the circle
and under the action of the natural involution.  Let \(\m{C}^+
\subseteq \m{C}\) denote the sub-cone consisting of positive definite
forms.

Let \(\m{T}\) denote the cone of positive, rapidly decreasing
sequences \((a_p)\) such that \(a_p = a_{-p}\) for all \(p \in \Z\).
Let \(\m{T}^+\) denote the sub-cone of strictly positive sequences.
\end{defn}

\begin{theorem}
\label{th:ipdesc}
The map \(\ipvc \to \left( \ipv{e_p}{e_p} \right)\) defines a
bijection of cones from \(\m{C}\) to \(\m{T}\) such that \(\m{C}^+\)
is carried onto \(\m{T}^+\).
\end{theorem}

\begin{proof}
As the set \(\{e_p : p \in \Z\}\) is a basis for \(\m{S}^*\), any
sesquilinear form, \(\ipvc\), on \(\m{S}^*\) is completely determined
by the \(\Z \times \Z\)-indexed set of numbers \(\{\ipv{e_p}{e_q}\}\).
We shall refer to this as the double sequence associated to
\(\ipvc\).

Suppose that \(\ipvc\) is a sesquilinear form on \(\m{S}^*\) invariant
under the action of \(R_\lambda\) for some \(\lambda \in S^1\) not of
finite order.  Then for all \(p,q \in \Z\), \(\ipv{R_\lambda
e_p}{R_\lambda e_q} = \ipv{e_p}{e_q}\).  Using the formula from
lemma~\ref{lem:opprop}, the left-hand side of this equation is
\(\lambda^{q-p} \ipv{e_p}{e_q}\).  As \(\lambda\) is not of finite
order, this implies that \(\ipv{e_p}{e_q} = 0\) for \(p \ne q\).  Thus
the double sequence associated to \(\ipvc\) is zero off the main
diagonal.

Conversely, suppose that \(\ipvc\) is a sesquilinear form on
\(\m{S}^*\) such that the associated double sequence is zero off the
main diagonal.  For \(a = (a^p) \in \m{S}^*\), the number
\(\ipv{a}{a}\) is given by the formula \(\sum
\abs{a^p}\ipv{e_p}{e_p}\).  Thus as \(R_\lambda a =
(\lambda^{-p}a^p)\), \(\ipv{R_\lambda a}{R_\lambda a} = \ipv{a}{a}\)
for any \(\lambda \in S^1\).  Hence \(\ipvc\) is invariant under the
circle action.

If, in addition, the natural involution acts unitarily -- that is, the
sesquilinear form is invariant under the action of the natural
involution -- then lemma~\ref{lem:opprop} shows that \(\ipv{e_p}{e_p}
= \ipv{e_{-p}}{e_{-p}}\).  The converse is immediate.

Let \(\ipvc\) be a sesquilinear form which is invariant under the
circle action and under the natural involution.  Let \(a_p =
\ipv{e_p}{e_p}\) for \(p \in \Z\).  The form \(\ipvc\) is continuous
and therefore defines a conjugate linear map \(\m{S}^* \to {\m{S}^*}^* =
\m{S}\).  Under this map, an element \(b \in \m{S}^*\) is taken to the
sequence \((\ipv{e_p}{b})\).  Let \(1\!\!1 \in \m{S}^*\) denote the
sequence consisting completely of \(1\)s.  Under the map \(\m{S}^* \to
\m{S}\) defined by the form, this element is taken to
\((\ipv{e_p}{1\!\!1}) = (a_p)\).  Hence the sequence \((a_p)\) is
rapidly decreasing and thus the map in the statement of the theorem is
well-defined.

Thus a sesquilinear form which is invariant under the circle action
and under the natural involution is completely determined by the
sequence \((\ipv{e_p}{e_p})\).  It is simple to see that the sequence
is positive if and only if the sesquilinear form is positive
semi-definite, and that the sequence is strictly positive if and only
if the sesquilinear form is positive definite.  Thus the sequence is
an element of \(\m{T}\) and is in \(\m{T}^+\) if and only if the
original sesquilinear form were positive definite.  Whence the map
\(\m{C} \to \m{T}\) is well-defined and injective.  A simple check
shows that this is a map of cones.

To show that the map is surjective, and hence a bijection, let \((a_p)
\in \m{T}\).  Let \(b = (b^p)\) and \(c = (c^p)\) be elements of
\(\m{S}^*\).  There exist integers \(m,n > 0\) such that
\((p^{-m}b^p)\) and \((p^{-n}c^p)\) are bounded.  As \((a_p)\) is
rapidly decreasing, the sequence \((p^{n+m+2}a_p)\) is bounded and
hence \((p^{n+m}a_p)\) is summable.  Hence \((b^p\conj{c^p}a_p)\) is a
summable sequence and thus the formula:
\[
(b,c) \to \sum_{p \in \Z} b^p \conj{c^p}a_p
\]
is well-defined as a sesquilinear map \(\m{S}^* \times \m{S}^* \to \C\).
It is evidently positive semi-definite.  To show continuity, it is
sufficient to show that it is continuous when restricted to each space
\(\{(x^p) : (p^{-n}x^p) \text{ is bounded}\}\).  Continuity of this
restriction follows from the estimate:
\[
\abs{ \sum_{p \in \Z} b^p \conj{c^p}a_p } \le \sup \{
\abs{p^{-n}b^p}\} \sup \{ \abs{ p^{-n} c^p} \} \sum_{p \in \Z}
\abs{p^{2n} a_p }.
\]

Thus the sequence \((a_p)\) defines a sesquilinear form on
\(\m{S}^*\).  It is clear that the associated double sequence for this
form is zero off the main diagonal and on the main diagonal is
\((a_p)\).  Thus it is invariant under the circle action and the
natural involution and so is an element of \(\m{C}\).  It is the
preimage of \((a_p)\) under the map \(\m{C} \to \m{T}\) showing that
the map is a bijection.
\end{proof}

Any continuous inner product on \(\m{S}^*\) defines a Hilbert space
completion, but the map from inner products to Hilbert space
completions is not injective.  Two inner products define the same
Hilbert space completion if and only if the identity map on
\(\m{S}^*\) extends to an isomorphism between the completions.  This
condition can be stated elegantly in terms of the  sequences in
\(\m{T}^+\) associated to the given inner products:

\begin{lemma}
  Let \((a_p), (b_p) \in \m{T}^+\).  The Hilbert space completions
  defined by the inner products associated to \((a_p)\) and \((b_p)\)
  are equivalent if and only if the sequences \((a_p/b_p)\) and
  \((b_p/a_p)\) are bounded.
\end{lemma}

We now turn to the operator \(z\) and determine the answer to the
following question: for which inner products on \(\m{S}^*\) is the
operator \(z\) continuous with respect to the inner product topology?

\begin{proposition}
  Let \((a_p) \in \m{T}^+\).  Let \(\ipvc\) be the associated inner
  product on \(\m{S}^*\).  The operator \(z\) is continuous with
  respect to the inner product topology if and only if the sequence of
  ratios \((a_p/a_{p+1})\) is bounded.

  In this case, \(\norm[z]^2 = \sup\{a_p/a_{p+1}\}\).
\end{proposition}

Notice that as \(a_p = a_{-p}\), the sequence \((a_p/a_{p-1})\) is
just \((a_p/a_{p+1})\) in reverse order.  Moreover, we cannot have
\(z\) acting unitarily as this would imply that \((a_p)\) is constant,
contradicting the fact that it is rapidly decreasing.

\begin{proof}
  Let \(\norm\) be the norm defined by the inner product.  Suppose
  that \(z\) is continuous with respect to the inner product topology
  on \(\m{S}^*\).  In particular, \(\norm[z e_{p+1}] \le
  \norm[z]\norm[e_{p+1}]\) for all \(p\).  From
  lemma~\ref{lem:opprop}, \(z e_{p+1} = e_p\).  Thus for \(p \in \Z\),
  \(\sqrt{a_p} \le \norm[z] \sqrt{a_{p+1}}\).  Hence the sequence
  \((a_p/a_{p+1})\) is bounded above by \(\norm[z]^2\).

  Conversely, suppose that \((a_p/a_{p+1})\) is bounded above by, say,
  \(M\).  Let \(b = (b^p) \in \m{S}^*\), then:
  \[
  \norm[z b]^2 = \sum \abs{b^{p+1}} a_p \le \sum \abs{b^{p+1}} Ma_{p+1}
  = M \norm[b]^2.
  \]
  Thus \(z\) is continuous with respect to \(\norm\) and so extends to
  a continuous linear operator on \(H\).  Moreover, \(\norm[z]^2 \le
  M\).

  Combining the two relationships for \(\norm[z]\) shows that
  \(\norm[z]^2 = \sup \{a_p/a_{p+1} : p \in \Z\}\) when either side
  exists.
\end{proof}

\begin{corollary}
  \label{cor:zq}
  Let \((a_p) \in \m{T}\) be such that \((a_p/a_{p+1})\) is bounded.
  For each \(q \in \Z\), the operator \(z^q\) is continuous with
  respect to the inner product topology and \(\norm[z^q]^2 =
  \sup\{a_p/a_{p+q}\}\).
\end{corollary}

\begin{corollary}
\label{cor:cvxip}
  Let \(\m{C}_z\) be the subset of \(\m{C}\) consisting of those inner
  products for which the operation of multiplication by \(z\) is
  continuous, \(\m{T}_z\) the corresponding sub-cone of \(\m{T}\).
  Then \(\m{C}_z\) is a non-empty sub-cone of \(\m{C}^+\).
\end{corollary}

\begin{proof}
  The set \(\m{T}_z\) consists of those sequences \((a_p) \in
  \m{T}^+\) for which \((a_p/a_{p+1})\) is bounded.  This is non-empty
  as the sequence \((2^{-\abs{p}})\) lies in \(\m{T}_z\).

  Clearly, if \((a_p) \in \m{T}_z\) then for any \(t > 0\), \((t a_p)
  \in \m{T}_z\).  If \((a_p), (b_p) \in \m{T}_z\) then there exist
  \(M,N > 0\) such that \(a_p/a_{p+1} \le M\) and \(b_p/b_{p+1} \le
  N\) for all \(p\).  Equivalently, \(a_p \le M a_{p+1}\) and \(b_p
  \le N b_{p+1}\).  Let \(R = \max\{M,N\}\), then \(a_p + b_p \le
  R(a_{p+1} + b_{p+1})\) so \(\left((a_p + b_p)/(a_{p+1} +
  b_{p+1})\right)\) is bounded, hence lies in \(\m{T}_z\).

  Therefore \(\m{T}_z\) is a sub-cone of \(\m{T}^+\) and so \(\m{C}_z\)
  is a sub-cone of \(\m{C}^+\).
\end{proof}

These inner products transfer to \(L^*\C\) via the isomorphism \(L^*\C
\to \m{S}^*\).  We can find a formula which is more natural on
\(L\C\).

\begin{proposition}
  Let \(\ipvc \in \m{C}\).  Let \((a_p) \in \m{T}\) be the associated
  sequence.  Thinking of \(\m{T}\) as a subset of \(\m{S}\), let
  \(\gamma_a \in L\C\) be the image of \((a_p)\) under the isomorphism
  \(\m{S} \cong L\C\).

  Under the isomorphism \(\m{S}^* \cong L^*\C\), the form \(\ipvc\) is
  given by the formula \((b,c) \to b ( \conj{c} \diamond \gamma_a)\)
  where \(c \diamond \gamma_a \in L\C\) is the map \(\lambda \to
  c(R_{\lambda^{-1}} \gamma_a)\).
\end{proposition}

\begin{proof}
  The map \(S^1 \times S^1 \to \C\) defined by \((\lambda, \mu) \to
  \gamma_a(\lambda^{-1}\mu)\) is the composition of smooth maps hence
  is smooth.  Therefore by the exponential law for smooth maps,
  \cite[I.3]{akpm}, its adjoint, \(\lambda \to
  R_{\lambda^{-1}}\gamma_a\), is a smooth map \(S^1 \to L\C\).  The
  element \(c \in L^*\C\) is a continuous linear map \(L\C \to \C\),
  hence is smooth, so the map \(\lambda \to
  c(R_{\lambda^{-1}}\gamma_a)\) is a smooth map \(S^1 \to \C\).  Thus
  the formula \((b,c) \to b(\conj{c} \diamond \gamma_a)\) makes sense.
  It is also evident that the map \(c \to c \diamond \gamma_a\) is
  continuous and so \((b,c) \to b( \conj{c} \diamond \gamma_a)\) is at
  least separately continuous and thus completely determined by its
  effect on a basis.

  Using the isomorphisms \(L\C \cong \m{S}\) and \(L^*\C \cong
  \m{S}^*\), we transfer these operators to \(\m{S}\) and \(\m{S}^*\).
  Under these isomorphisms, \(R_{\lambda^{-1}} \gamma_a\) becomes the
  sequence \((\lambda^{-p}a_p)\) and so \(e_q(R_{\lambda^{-1}}
  \gamma_a) = \lambda^{-q}a_q\).  Thus \(e_q \diamond \gamma_a\) is
  the sequence corresponding to the function \(\lambda \to
  \lambda^{-q}a_q\) which is \(a_q e^{-q}\).  Therefore, \(\Sigma
  e_p(e_q \diamond \gamma_a) = e_{-p} (a_q e^{-q}) = a_q \delta_p^q\).

  Hence the sesquilinear form on \(L^*\C\), \((b,c) \to b( \conj{c}
  \diamond \gamma_a)\), corresponds to the original sesquilinear form
  on \(\m{S}^*\), \(((b^p), (c^p)) \to \sum b^p \conj{c^p} a_p\).
\end{proof}

The inner products we consider on \(\m{S}^*\) and \(L^*\C\) arise as
inner products on the underlying real spaces and therefore give a
classification of inner products on \(L^*\R\) which are invariant under
the circle action and the natural involution, and also of those for
which the operations of multiplication by \(\cos \theta\) and by
\(\sin \theta\) are continuous.

\begin{proposition}
  The sesquilinear forms on \(\m{S}^*\) and \(L^*\C\) considered above
  are the complexifications of sesquilinear forms on the underlying
  real spaces of both \(\m{S}^*\) and \(L^*\C\).
\end{proposition}

\begin{proof}
  For \(\m{S}^*\), this is evident from the formula.  For \(L^*\C\), it
  follows from the invariance under \(\iota\) together with the fact
  that \(\iota\) intertwines the complex conjugation operators arising
  from \(\m{S}^*\) and \(L^*\C\) (note that the isomorphism \(L^* \C \to
  \m{S}^*\) does not induce an isomorphism of real structures and thus
  the complex conjugation operators differ).
\end{proof}

\subsection{Polarisations}

In this section we examine how the theory of polarisations, and thus
of unitary structures, interacts with these inner products on the
space of distributions.  We examine an inner product on \(L^*\C^n\)
determined by a sequence in \(\m{T}_z\).  To pass from an inner
product on \(L^*\C\) to one on \(L^*\C^n\), we use the isomorphism
\(L^*\C^n \cong L^*\C \otimes \C^n\) together with the the standard
inner product on \(\C^n\).

\begin{lemma}
  \label{lem:polop}
  Let \(J\) be the operator on \(\m{S}^*\) defined by \(J e_p = -
  (-1)^{\text{sign}(p)} i e_p\).  Let \(\ipvc \in \m{C}^+\) be an
  inner product on \(\m{S}^*\) and let \(H\) be the corresponding
  Hilbert space completion.  Then the operator \(J\) defines a
  polarisation of \(H\).

  This extends in a natural way to a polarisation of the Hilbert
  completion of \(L^* \C^n\).
\end{lemma}

\begin{proof}
  Let \((a_p) \in \m{T}^+\) be the sequence corresponding to the inner
  product.  For \(b = (b^p) \in \m{S}^*\), it follows straight from
  the formula for \(J\) that \(\ipv{J b}{J b} = \ipv{b}{b}\) and
  therefore \(J\) extends to a unitary operator on \(H\).  It
  satisfies \(J^2 = -1\) and \(J \pm i I\) are not finite rank.
  Therefore, it defines a polarisation on \(H\).

  To extend this to the Hilbert completion of \(L^* \C^n\), we
  observe that this completion is naturally isomorphic to \(H \otimes
  \C^n\).  The polarising operator \(J\) on \(H\) defines one on \(H
  \otimes \C^n\) by taking \(J \otimes I_n\).
\end{proof}

\begin{defn}
  The polarisation \(\m{J}\) so defined on the completion of \(L^*
  \C^n\) is called the \emph{standard polarisation}.
\end{defn}

\begin{proposition}
  Let \((a_p) \in \m{T}_z\).  Let \(H\) be the associated Hilbert
  space completion of \(L^* \C^n\).  The polynomial loop group
  \(L_\pol U_n\) acts continuously on \(H\) and preserves the
  polarisation.
\end{proposition}

\begin{proof}
  Let \(J\) be the polarising operator on \(H\) as defined in lemma
  \ref{lem:polop}.  Let \(\m{L}_\m{J}(H)\) be the set of all bounded
  linear operators \(A\) on \(H\) such that \([A,J]\) is
  Hilbert-Schmidt.  It is clear that \(\gl_\m{J}(H) = \gl(H) \cap
  \m{L}_\m{J}(H)\).  The norm of an element \(A \in \m{L}_\m{J}(H)\)
  is the sum of the operator norm of of \(A\) and the Hilbert-Schmidt
  norm of \([A,J]\).

  There is an isometry \(M_n(\C) \to \m{L}(H)\) given by \(A(a \otimes
  v) = a \otimes Av\), thinking of \(H\) as the completion of \(L^*
  \C \otimes \C^n\).  This maps continuously into
  \(\m{L}_\m{J}(H)\) since \(A \in M_n(\C)\) commutes with \(J\).

  The operator \(z\) acts continuously on \(H\) and \([J,z]\) is
  finite rank.  It therefore lies in \(\m{L}_\m{J}(H)\).  Thus
  the image of \(L_\pol M_n(\C)\) in \(\m{L}(H)\) lies in
  \(\m{L}_\m{J}(H)\).  Thus the image of \(L_\pol U_n\) lies in
  \(\gl_\m{J}(H)\).
\end{proof}

\begin{proposition}
  The inclusion \(L_\pol U_n \to \gl_\m{J}(H)\) is homotopic to the
  standard inclusion which factors through \(L U_n\).
\end{proposition}

\begin{proof}
  Let \(T : H \to L^{2,*} \C^n\) be the isometry which takes \(e_p\)
  to \(\sqrt{a_p} e_p\).  This identifies \(\gl_\m{J}(H)\) with
  \(\gl_\m{J}(L^{2,*} \C^n)\) and so defines the map \(L_\pol U_n \to
  \gl_\m{J}(L^{2,*} \C^n)\).

  Let \(\zeta_t : L^{2,*} \C^n \to L^{2,*} \C^n\) be the map defined
  by \(\zeta_t(e_p) = (a_{p-1}/a_p)^{t/2} e_{p-1}\).  As \(a_p\) is
  positive for all \(p\), the formula makes sense.  Since \((a_p)\)
  lies in \(\m{T}_z\), \((a_{p-1}/a_p)\) is bounded above and below so
  \(\zeta_t\) is an isomorphism of Hilbert spaces.

  The map \(\zeta_0\) is the (adjoint of the) map \(z\).  The map
  \(\zeta_1\) is the map \(T z T^{-1}\).  Therefore the two inclusions
  of \(L_\pol U_n\) in \(\gl_\m{J}(H)\) are \(\sum z^q A_q \to \sum
  \zeta_0^q A_q\) and \(\sum z^q A_q \to \sum \zeta_1^q A_q\).  The
  required homotopy is \(F(\sum z^q A_q, t) = \sum \zeta_t^q A_q\).
\end{proof}

This homotopy equivalence is closely related to the deformation
retract \(\gl(H) \to U(H)\).  This retract uses the polar
decomposition: any invertible linear operator on \(H\) can be written
in the form \(A = Q \abs{A}\) where \(\abs{A}\) is a self-adjoint
operator with strictly positive eigenvalues and \(Q\) is a unitary
operator.  The retraction maps \(A\) to \(Q\) and the homotopy
equivalence is given by contracting the eigenvalues of \(\abs{A}\) to
\(1\).

This retract and homotopy equivalence restricts to give retracts and
homotopy equivalences of the various standard subgroups of \(\gl(H)\)
onto their unitary counterparts.  Writing \(H_\R\) for an underlying
real Hilbert space of \(H\), the map \(\gl(H) \to U(H)\) also maps:
\begin{align*}
\gl(H_\R) &\to O(H_\R), \\
\gl_\m{J}(H) &\to U_\m{J}(H), \\
\gl_J(H_\R) &\to O_J (H_\R), \\
\gl_n(\C) &\to U_n, \\
\gl_n(\R) &\to O_n.
\end{align*}
The relation with the homotopy equivalence above comes about because
in the polar decomposition of \(\zeta_1\), the unitary operator is
\(z\).  The homotopy \(\zeta_t\) is the homotopy which contracts the
positive part to the identity.

\label{page:last}

\end{document}